\documentclass[10.5pt]{article}
\usepackage[leqno]{amsmath}
\usepackage{amsfonts}
\usepackage{graphicx}

\usepackage{amsmath}
\usepackage{amssymb}
\usepackage{latexsym}
\usepackage{amsmath, amsfonts,amssymb, amsthm, euscript,makeidx,color,mathrsfs}

\def\sqr#1#2{{\vcenter{\vbox{\hrule height.#2pt
              \hbox{\vrule width.#2pt height#1pt \kern#1pt \vrule width.#2pt}
              \hrule height.#2pt}}}}
\def\signed #1{{\unskip\nobreak\hfil\penalty50
              \hskip2em\hbox{}\nobreak\hfil#1
              \parfillskip=0pt \finalhyphendemerits=0 \par}}
\def\endpf{\signed {$\sqr69$}}

\def\3n{\negthinspace \negthinspace \negthinspace }
\def\2n{\negthinspace \negthinspace }
\def\1n{\negthinspace }

\def\dbI{\mathbb{I}}

\def\dbR{\mathbb{R}}
\def\dbS{\mathbb{S}}

\def\={\buildrel \triangle \over =}

\def\ds{\displaystyle}

\def\ns{\noalign{\ss}}
\def\a{\alpha}
\def\b{\beta}

\def\d{\delta}
\def\e{\varepsilon}

\def\l{\lambda}
\def\m{\mu}
\def\n{\nu}

\def\t{\tau}

\def\o{\omega}

\def\F{\Phi}

\def\cB{{\cal B}}

\def\cU{{\cal U}}

\def\cX{{\cal X}}

\def\BA{{\bf A}}

\def\BC{{\bf C}}

\def\BF{{\bf F}}

\def\ss{\smallskip}
\def\ms{\medskip}

\def\q{\quad}
\def\qq{\qquad}
\def\hb{\hbox}

\def\lan{\mathop{\langle}}
\def\ran{\mathop{\rangle}}

\def\h{\widehat}
\def\wt{\widetilde}

\def\cd{\cdot}
\def\cds{\cdots}

\def\ae{\hbox{\rm a.e.{ }}}

\def\({\Big (}
\def\){\Big )}
\def\[{\Big[}
\def\]{\Big]}
\def\bde{\begin{definition}}
\def\ede{\end{definition}}
\def\be{\begin{equation}}
\def\bel{\begin{equation}\label}
\def\ee{\end{equation}}
\def\bt{\begin{theorem}}
\def\et{\end{theorem}}
\def\bc{\begin{corollary}}
\def\ec{\end{corollary}}
\def\bl{\begin{lemma}}
\def\el{\end{lemma}}
\def\bp{\begin{proposition}}
\def\ep{\end{proposition}}
\def\bas{\begin{assumption}}
\def\eas{\end{assumption}}
\def\br{\begin{remark}}
\def\er{\end{remark}}
\def\ba{\begin{array}}
\def\ea{\end{array}}
\def\ed{\end{document}}

\def\square#1{\vbox{\hrule\hbox{\vrule height#1%
     \kern#1\vrule}\hrule}}
\def\rectangle#1#2{\vbox{\hrule\hbox{\vrule height#1%
     \kern#2\vrule}\hrule}}

\font\tenbb=msbm10 \font\sevenbb=msbm7 \font\fivebb=msbm5

\newfam\bbfam
\scriptscriptfont\bbfam=\fivebb \textfont\bbfam=\tenbb
\scriptfont\bbfam=\sevenbb

\newtheorem{lemma}{Lemma}[section]
\newtheorem{remark}{Remark}[section]

\newtheorem{theorem}{Theorem}[section]
\newtheorem{corollary}{Corollary}[section]

\newtheorem{definition}{Definition}[section]
\newtheorem{proposition}{Proposition}[section]
\newtheorem{assumption}{Assumption}[section]

\makeatletter
   
   \@addtoreset{equation}{section}
\makeatother

\begin{document}
\title{\bf A Deterministic Affine-Quadratic\\ Optimal Control Problem\footnote{
The first author was supported in part by NSFC under Grant 71163046
and China State Scholarship Fund under Grant [2009]5004, the second
author was supported in part by NSF under Grant DMS-1007514.}}

\author{Yuanchang Wang$^{a,b}$~~and~~Jiongmin Yong$^b$ \\
$^a$ School of Mathematics, Yunnan Normal University , Kunming,
650500, China.\\
$^b$ Department of Mathematics, University of Central Florida,
Orlando, FL 32816, USA.}

\maketitle

\begin{abstract}
A Deterministic affine quadratic optimal control problem is
considered. Due to the nature of the problem, optimal controls exist
under some very mild conditions. Further, it is shown that under
some assumptions, the value function is differentiable and therefore
satisfies the corresponding Hamilton-Jacobi-Bellman equation in the
classical sense. Moreover, the so-called quasi-Riccati equation is
derived and any optimal control admits a state feedback
representation.
\end{abstract}
\ms


\section{Introduction.}\label{1}

Consider the following controlled ordinary differential equation
(ODE, for short):
\bel{1.1}\left\{\2n\ba{ll}
\ns\ds\dot X(s)=A(s,X(s))+B(s,X(s))u(s),\qq s\in[t,T],\\
\ns\ds X(t)=x,\ea\right.\ee
with cost functional
\bel{1.2}\ba{ll}
\ns\ds J(t,x;u(\cd))\1n=\2n\int_t^T\2n\big[Q(s,X(s))\1n+\1n\lan
S(s,X(s)),u(s)\ran+{1\over2}\lan
R(s,X(s))u(s),u(s)\ran\big]ds\1n+\1n G(X(T)),\ea\ee
where $A:[0,T]\times\dbR^n\to\dbR^n$,
$B:[0,T]\times\dbR^n\to\dbR^{n\times m}$,
$Q:[0,T]\times\dbR^n\to\dbR$, $S:[0,T]\times\dbR^n\to\dbR^m$,
$R:[0,T]\times\dbR^n\to\dbS_+^m$ ($\dbS^m$ is the set of all
symmetric matrices, and $\dbS^m_+$ is the set of all positive
semi-definite matrices), and $G:[0,T]\times\dbR^n\to\dbR$ are some
given maps. Let $\cU[t,T]$ be the set of all admissible controls
(which will be specified in the next section) on $[t,T]$. Under some
mild conditions, for any $(t,x)\in[0,T)\times\dbR^n$ and
$u(\cd)\in\cU[t,T]$, the state equation (\ref{1.1}) admits a unique
solution $X(\cd)\equiv X(\cd\,;t,x,u(\cd))$ and the cost functional
(\ref{1.2}) is well-defined. Then we can pose the following optimal
control problem.

\ms

\bf Problem (AQ). \rm For any given $(t,x)\in[0,T)\times\dbR^n$,
find a $u^*(\cd)\in\cU[t,T]$ such that
\bel{}J(t,x;u^*(\cd))=\inf_{u(\cd)\in\cU[t,T]}J(t,x;u(\cd))\equiv
V(t,x).\ee
Any $u^*(\cd)$ satisfying the above is called an {\it optimal
control} for $(t,x)$, and the corresponding $X^*(\cd)\equiv
X(\cd\,;t,x,u^*(\cd))$ is called an {\it optimal trajectory} for
$(t,x)$. The pair $(X^*(\cd),u^*(\cd))$ is called an {\it optimal
pair} of Problem (AQ) for the initial pair $(t,x)$. The function
$V(\cd\,,\cd)$ is called the {\it value function} of Problem (AQ).

\ms

We note that the right hand side of the state equation is {\it
affine} with respect to the control and the integrand in the cost
functional is up to quadratic with respect to the control.
Therefore, we call such a problem an {\it affine-quadratic} optimal
control problem (AQ problem, for short). We see that if
\bel{}\left\{\2n\ba{ll}
\ns\ds A(t,x)=A(t)x,\q B(t,x)=B(t),\q Q(t,x)={1\over2}\lan Q(t)x,x\ran,\\
\ns\ds S(t,x)=S(t)x,\q R(t,x)=R(t),\q G(x)={1\over2}\lan
Gx,x\ran,\ea\right.\qq\forall(t,x)\in[0,T]\times\dbR^n,\ee
for some matrix-valued functions $A(\cd)$, $B(\cd)$, $Q(\cd)$,
$S(\cd)$, $R(\cd)$, and some matrix $G$, then our Problem (AQ) is
reduced to a standard linear-quadratic optimal control problem (LQ
problem, for short).

\ms

It is well-known that for LQ problem, under suitable conditions, one
has the existence of a unique optimal control which admits a state
feedback representation via the solution of a differential Riccati
equation (\cite{Kalman 1960}, see also \cite{Yong-Zhou 1999}). On
the other hand, for optimal control problem of general nonlinear
ordinary differential equation with a Bolza type cost functional,
one generally does not expect the existence of an optimal control;
However, under some mild conditions, one can characterize the value
function of the optimal control problem as the unique viscosity
solution to the so-called Hamilton-Jacobi-Bellman (HJB, for short)
equation (\cite{Bardi--Capuzzo-Dolcetta 1997}, see also
\cite{Bardi-Da Lio 1997}, \cite{Qiu-Yong 2012}, and the references
cited therein). Note that our Problem (AQ) is between general
(nonlinear) optimal control problems and LQ problems. Therefore, one
expects some results ``between'' the results for the above-mentioned
two kinds of problems. A little more precisely, under certain
conditions, we will have the existence of optimal controls. Further,
it is possible to have state feedback representation of optimal
control via a solution to the so-called quasi-Riccati equation. We
would like to mention that Problem (AQ) with the state equation
being linear and with the maps $x\mapsto Q(t,x)$ and $x\mapsto G(x)$
being convex, and $S(t,x)\equiv0$ was studied in \cite{You 1987} and
\cite{You 1997} by means of the quasi-Riccati equations. Also,
without giving details, Problem (AQ) for stochastic differential
equations was briefly discussed in \cite{Yong 2002}.

\ms

Our approach is a combination of variational method and dynamic
programming method. The key is to obtain, under certain hypotheses,
the convexity of the map $u(\cd)\mapsto J(t,x;u(\cd))$ which will
lead to the differentiability of the value function $V(t,x)$. Then
the Hamilton-Jacobi-Bellman (HJB, for short) equation will be
satisfied in the classical sense. Furthermore, by differentiating
the HJB equation we obtain a quasi-Riccati equation.

\ms

We refer to \cite{Frankowska} and \cite{Bardi--Capuzzo-Dolcetta
1997} for excellent surveys on the value function of optimal control
theory. See also \cite{Benveniste-Scheinkman
1979,Cannarsa-Frankowska 1991,Bardi-Da Lio
1997,Rincon-Zapatero--Santos} for some relevant results concerning
the differentiability of value functions.

\ms

The rest of the paper is organized as follows. Section 2 collects
some preliminary results. In Section 3, we present the existence of
optimal controls for our Problem (AQ) and recall a Pontryagin type
minimum principle. In Section 4, we derive the first and the second
order variations of the cost functional with respective to the
control. The invertibility of the Hessian $D_{uu}J(t,x;u(\cd))$ of
the cost functional with respect to the control variable is obtained
in Section 5, under certain sufficient conditions. In Section 6, we
derive the so-called quasi-Riccati equation in a very natural way,
via which a state feedback representation of the optimal control is
obtained. A couple of illustrative examples are presented as well.
Finally, some concluding remarks are collected in Section 7.

\section{Preliminaries.}

Throughout this paper, we let $U\subseteq\dbR^m$ be a nonempty
convex and closed set, not necessarily bounded. For convenience, we
assume hereafter that $0\in U$. Note that it could be $U=\dbR^m$.
Now, we introduce the following standing assumptions.

\ms

{\bf(H1)} The maps $A:[0,T]\times\dbR^n\to\dbR^n$ and
$B:[0,T]\times\dbR^n\to\dbR^{n\times m}$ are continuous. There exist
constants $L_A,L_B,\wt L_B>0$ such that
\bel{2.1}|A(t,x)-A(t,\bar x)|\le L_A|x-\bar x|,\qq\forall
t\in[0,T],~x,\bar x\in\dbR^n,\ee
\bel{2.2}|B(t,x)-B(t,\bar x)|\le\wt L_B|x-\bar x|,\qq\forall
t\in[0,T],~x,\bar x\in\dbR^n,\ee
and
\bel{2.3}\lan[B(t,x)-B(t,\bar x)]^T(x-\bar x),u\ran\le L_B|x-\bar
x|^2,\qq\forall(t,u)\in[0,T]\times U,~x,\bar x\in\dbR^n.\ee

\ms

Note that condition (\ref{2.3}) is equivalent to the following:
\bel{2.4}\sup_{u\in U,\,x,\bar x\in\dbR^n,\,x\ne\bar
x}{\lan[B(t,x)-B(t,\bar x)]^T(x-\bar x),u\ran\over|x-\bar x|^2}\le
L_B.\ee
On the other hand, under (\ref{2.2}), the set
\bel{}\cX=\Big\{{[B(t,x)-B(t,\bar x)]^T(x-\bar x)\over|x-\bar
x|^2}\Bigm|x,\bar x\in\dbR^n,~x\ne\bar x\Big\}\subseteq\cB^m_{\wt
L_B}(0),\ee
where $\cB^m_r(0)$ is the ball in $\dbR^m$ centered at $0$ with
radius $r$. Therefore, in the case $U$ is bounded, (\ref{2.4}) is
satisfied with
$$L_B\ge\wt L_B\sup_{u\in U}|u|.$$
In the case $U=\dbR^m$, (\ref{2.4}) is equivalent to the following:
\bel{2.10}[B(t,x)-B(t,\bar x)]^T(x-\bar x)=0,\qq\forall
t\in[0,T],~x,\bar x\in\dbR^n.\ee
If we denote
$$B(t,x)=\(B^1(t,x),B^2(t,x),\cds,B^m(t,x)\),\q B^i:[0,T]\times\dbR^n\to
\dbR^n,~1\le i\le m,$$
then (\ref{2.10}) is equivalent to the following:
$$\lan B^i(t,x)-B^i(t,\bar x),x-\bar x\ran=0,\qq1\le i\le m.$$
This is the case if $B_x^i(t,x)$ is skew symmetric, for each $1\le
i\le m$. In particular, this is the case, of course, if $B(t,x)$ is
independent of $x$. Note that even if $B(t,x)$ is independent of
$x$, due to the fact that $x\mapsto A(t,x)$ is not necessarily
linear, we still have a nonlinear state equation.

\ms

Next, we introduce the following hypothesis for the functions
appearing in the cost functional.

\ms

{\bf(H2)} Maps $Q:[0,T]\times\dbR^n\to\dbR$,
$S:[0,T]\times\dbR^n\to\dbR^m$, $R:[0,T]\times\dbR^n\to\dbS^m$, and
$G:\dbR^n\to\dbR$ are continuous. There are constants
$L,Q_0,G_0,S_0>0$, $\e_0\in(0,1)$, and a continuous function
$\rho:\dbR^n\to[\rho_0,\infty)$ with $\rho_0>0$ such that
\bel{}\left\{\2n\ba{ll}
\ns\ds R(t,x)\ge\rho(x)I,\\
\ns\ds(1-\e_0)Q(t,x)-{1\over2}S(t,x)^TR(t,x)^{-1}S(t,x),~G(x)\ge-L,\\
\ns\ds Q(t,x)\le Q_0(1+|x|^2),\q G(x)\le G_0(1+|x|^2),\q|S(t,x)|\le
S_0(1+|x|),\\
\ns\ds\qq\qq\qq\qq\qq\qq\qq\forall(t,x)\in[0,T]\times\dbR^n.\ea\right.\ee
\bel{2.8}\left\{\2n\ba{ll}
\ns\ds|A(t,x)+B(t,x)R(t,x)^{-1}S(t,x)|\le L(1+|x|),\\
\ns\ds|B(t,x)R(t,x)^{-1}B(t,x)^T|\le
L,\ea\right.\qq(t,x)\in[0,T]\times\dbR^n.\ee

\ms

We also need the following assumption later.

\ms

\bf(H3) \rm The map
$(t,x)\mapsto(A(t,x),B(t,x),Q(t,x),S(t,x),R(t,x),G(x))$ is twice
continuously differentiable.

\ms

For any $0\le t<T$, let
\bel{}\cU[t,T]=\Big\{u(\cd)\in L^2(t,T;\dbR^m)\bigm|u(s)\in U,\q\ae
s\in[t,T]\Big\}.\ee
Any $u(\cd)\in\cU[t,T]$ is called an admissible control on $[t,T]$.
We denote
$$\|u(\cd)\|_{L^2(t,s)}=\(\int_t^s|u(r)|^2dr\)^{1\over2},\qq\forall
u(\cd)\in\cU[t,s].$$

\ms

The following simple result is concerned with the well-posedness of
the state equation (\ref{1.1}), whose proof is straightforward.

\ms

\bf Proposition 2.1. \sl Let {\rm(H1)} hold. Then for any initial
pair $(t,x)\in[0,T]\times\dbR^n$ and $u(\cd)\in\cU[t,T]$, equation
$(\ref{1.1})$ admits a unique solution $X(\cd)\equiv
X(\cd\,;t,x,u(\cd))$, and the following estimate holds:
\bel{2.7a}|X(s;t,x,u(\cd))|\le
K\[1+|x|+\|u(\cd)\|_{L^2(t,s)}\],\q\forall s\in[t,T],\ee
and
\bel{2.7b}|X(s;t,x,u(\cd))-x|\le
K\[1+|x|+\|u(\cd)\|_{L^2(t,s)}\]\[\sqrt{s-t}+\|u(\cd)\|_{L^2(t,s)}\]\sqrt{s-t},\ee
hereafter, $K>0$ denotes a generic constant which can be different
from line to line. Further, for any $t\in[0,T]$, $x,\bar
x\in\dbR^n$, and $u(\cd)\in\cU[t,T]$, it holds
\bel{2.17}|X(s;t,x,u(\cd))-X(s;t,\bar x,u(\cd))|\le e^{(L_A+
L_B)(T-t)}|x-\bar x|,\qq s\in[t,T].\ee

\ms

\rm

As a consequence of the above, using the technique found in
\cite{Qiu-Yong 2012}, we have the following result on the value
function.

\ms

\bf Proposition 2.2. \sl Let {\rm(H1)}--{\rm(H2)} hold. Then the
value function $V(\cd\,,\cd)$ is continuous and there exists a
constant $K>0$ such that
\bel{}-L(T-t+1)\le V(t,x)\le
K(1+|x|^2),\qq\forall(t,x)\in[0,T]\times\dbR^n,\ee
and
\bel{}|V(t,x)-V(t,\bar x)|\le K(|x|\vee|\bar x|)|x-\bar
x|,\qq\forall t\in[0,T],~x,\bar x\in\dbR^n,\ee
where $|x|\vee|\bar x|=\max\{|x|,|\bar x|\}$. Moreover, the value
function $V(\cd\,,\cd)$ is the unique viscosity solution to the
following HJB equation:
\bel{}\left\{\2n\ba{ll}
\ns\ds V_t(t,x)+\lan V_x(t,x),A(t,x)\ran+Q(t,x)\\
\ns\ds\qq+\inf_{u\in U}\[\lan B(t,x)^TV_x(t,x)+S(t,x),u\ran
+{1\over2}\lan R(t,x)u,u\ran\]\1n=0,
\q(t,x)\in[0,T]\times\dbR^n,\\
\ns\ds V(T,x)=G(x).\ea\right.\ee

\ms

\rm

Note that in the case $U=\dbR^m$, the above HJB equation can be
written as
\bel{}\left\{\2n\ba{ll}
\ns\ds
V_t(t,x)+H(t,x,V_x(t,x))=0,\qq(t,x)\in[0,T]\times\dbR^n,\\
\ns\ds V(T,x)=G(x),\qq x\in\dbR^n,\ea\right.\ee
with
\bel{H}\ba{ll}
\ns\ds H(t,x,p)\1n=\1n Q(t,x)\1n+\1n\lan
p,A(t,x)\ran\1n-{1\over2}\(B(t,x)^T\1n p\1n+\1n S(t,x)\)^T\1n
R(t,x)^{-1}
\(B(t,x)^T\1n p\1n+\1n S(t,x)\),\\
\ns\ds\qq\qq\qq\qq\qq\qq\qq\qq(t,x,p)\in[0,T]\times\dbR^n\times\dbR^n.\ea\ee
Further, in the case that $V(t,x)$ is differentiable, it is the
classical solution to the above HJB equation and the optimal control
admits the following representation:
$$u(s)=-R(s,X(s))^{-1}\[B(s,X(s))^TV_x(s,X(s))+S(s,X(s))\],\qq s\in[t,T],$$
with $X(\cd)$ being the solution to the closed-loop system:
$$\left\{\2n\ba{ll}
\ns\ds\dot
X(s)=A(s,X(s))-B(s,X(s))R(s,X(s))^{-1}\[B(s,X(s))^TV_x(s,X(s))+S(s,X(s))\],
\q s\in[t,T],\\
\ns\ds X(t)=x.\ea\right.$$

\ms

\rm

From \cite{Qiu-Yong 2012}, we note that to guarantee the uniqueness
of viscosity solution to the HJB equation, we need
\bel{2.17}\left\{\2n\ba{ll}
\ns\ds|H(t,x,p)-H(t,y,p)|\le\o(|x|+|y|,|p|,|x-y|),\qq
t\in[0,T],x,y,p\in\dbR^n,\\
\ns\ds|H(t,x,p)-H(t,x,q)|\le K_0\sum_{i=1}^k\lan
x\ran{}^{\l_i}\big(|p|\vee|q|\big)^{\n_i}|p-q|,\qq
t\in[0,T],x,p,q\in\dbR^n,\\
\ns\ds|G(x)-G(y)|\le K_0\big(\lan x\ran\vee\lan
y\ran\big)^{\m-1}|x-y|,\qq\forall x,y\in\dbR^n,\\
\ns\ds\l_i,\n_i\ge0,\q\l_i+(\m-1)\n_i\le1,\qq1\le i\le
k,\ea\right.\ee
with $\lan x\ran=\sqrt{1+|x|^2}$. For the current case, we may let
$\m=2$. Then
$$|G(x)-G(y)|\le L(\lan x\ran\vee\lan y\ran)|x-y|,\qq\forall
x,y\in\dbR^n.$$
When $U=\dbR^m$, the Hamiltonian has the explicit form (\ref{H}).
Clearly, the first condition in (\ref{2.17}) holds. For the second
condition, we observe that
$$\ba{ll}
\ns\ds|H_p(t,x)|\le|A(t,x)+B(t,x)R(t,x)^{-1}S(t,x)|+|B(t,x)R(t,x)^{-1}B(t,x)^Tp|\\
\ns\ds\qq\qq\le K_0\big(\1n\lan x\ran+|p|\,\big),\ea$$
which is implied by (\ref{2.8}). Thus, the second condition holds
with
$$\l_1=\n_2=1,\qq\l_2=\n_1=0.$$

\section{Existence of Optimal Controls and Minimum Principle.}

We first present the following result.

\ms

\bf Proposition 3.1. \rm Under (H1)--(H2), for any initial pair
$(t,x)\in[0,T)\times\dbR^n$, Problem (AQ) admits an optimal control.

\ms

\it Proof. \rm Let $(t,x)\in[0,T)\times\dbR^n$ be given. Let
$X_0(\cd)=X(\cd\,;t,x,0)$. According to (\ref{2.7a}) , we have
\bel{}|X_0(s)|\le K(1+|x|),\qq\forall s\in[t,T],\ee
for some $K>0$. Let $u^k(\cd)\in\cU[t,T]$ be a minimizing sequence
with the corresponding state trajectory $X^k(\cd)\equiv
X(\cd\,;t,x,u^k(\cd))$. Then we may assume that
$$\ba{ll}
\ns\ds J(t,x;0)+1\ge J(t,x;u^k(\cd))\\
\ns\ds=\int_t^T\[Q(s,X^k(s))+\lan
S(s,X^k(s)),u^k(s)\ran+{1\over2}\lan
R(s,X^k(s))u^k(s),u^k(s)\ran\]ds+G(X^k(T))\\
\ns\ds=\int_t^T\[{1\over1-\e_0}\((1-\e_0)Q(s,X^k(s))-{1\over2}
S(s,X^k(s))^TR(s,X^k(s))^{-1}S(s,X^k(s))\)\\
\ns\ds\qq\qq+{1\over2}\big|(1-\e_0)^{1\over2}R(s,X^k(s))^{1\over2}u^k(s)
+(1-\e_0)^{-{1\over2}}R(s,X^k(s))^{-{1\over2}}S(s,X^k(s))\big|^2\\
\ns\ds\qq\qq+{\e_0\over2}\lan R(s,X^k(s))u^k(s),u^k(s)\ran\]ds+G(X^k(T))\\
\ns\ds\ge-{L(T-t)\over1-\e_0}+{\e_0\rho_0\over2}\int_t^T|u^k(s)|^2ds-L.\ea$$
Thus,
\bel{}\int_t^T|u^k(s)|^2ds\le K,\qq\forall k\ge1.\ee
Consequently,
$$|X^k(s)|\le K\((1+|x|+\|u^k(\cd)\|_{\cU[t,T]}\)\le K,\qq\forall
s\in[t,T],~k\ge1.$$
Then for any $t\le s<\t\le T$,
$$\ba{ll}
\ns\ds|X^k(\t)-X^k(s)|\le\int_s^\t\(A_0+L_A|X^k(r)|+(B_0+\wt L_B|X^k(r))|)
|u^k(r)|\)dr\\
\ns\ds\qq\qq\qq\q\;\le
K(\t-s)+K(\t-s)^{1\over2}\|u^k(\cd)\|_{\cU[t,T]}\le
K(\t-s)^{1\over2}.\ea$$
Thus, $\{X^k(\cd)\}$ is uniformly bounded and equicontinuous. Hence,
we may assume that $X^k(\cd)\to X^*(\cd)$ in $C([t,T];\dbR^n)$. Then
a standard argument applies to get the existence of an optimal
control (see \cite{Berkovitz}). \endpf

\ms

Now, we have the following necessary conditions for any optimal pair
of Problem (AQ).

\ms

\bf Proposition 3.2. \sl Let {\rm(H1)}--{\rm(H3)} hold and
$(t,x)\in[0,T)\times\dbR^n$ be given. Let $(X^*(\cd),u^*(\cd))$ be
an optimal pair of Problem (AQ) for $(t,x)$. Then the following
adjoint equation admits a unique solution
\bel{3.1}\left\{\2n\ba{ll}
\ns\ds\dot Y(s)=-\[A_x(s,X^*(s))+\sum_{j=1}^m
u^*_j(s)B^j_x(s,X^*(s))\]^TY(s)-Q_x(s,X^*(s))^T\\
\ns\ds\qq\qq-S_x(s,X^*(s))^Tu^*(s)-{1\over2}\sum_{j,k=1}^m
u^*_j(s)u^*_k(s)R^{jk}_x(s,X^*(s))^T,\qq s\in[t,T],\\
\ns\ds Y(T)=G_x(X^*(T))^T,\ea\right.\ee
and the following minimum condition holds:
\bel{}\ba{ll}
\ns\ds\[B(s,X^*(s))^TY(s)+S(s,X^*(s))\]u^*(s)+{1\over2}u^*(s)^TR(s,X^*(s))u^*(s)\\
\ns\ds=\min_{u\in
U}\Big\{\[B(s,X^*(s))^TY(s)+S(s,X^*(s))\]u+{1\over2}u^TR(s,X^*(s))u\Big\},\q
s\in[t,T].\ea\ee
In the above, $B(s,x)=(B^1(s,x),B^2(s,x),\cds,B^m(s,x))$ with
$B^i:[0,T]\times\dbR^n\to\dbR^n$, and
$B^i_x:[0,T]\times\dbR^n\to\dbR^{n\times n}$. In particular, if
$U=\dbR^m$, we have
\bel{3.5}u^*(s)=-R(s,X^*(s))^{-1}\[B(s,X^*(s))^TY(s)+S(s,X^*(s))\],\qq
s\in[t,T].\ee

\rm

\ms

From the above result, we see that under (H1)--(H3) with $U=\dbR^m$,
for any $(t,x)\in[0,T)\times\dbR^n$, the following coupled two-point
boundary value problem admits a solution $(X(\cd),Y(\cd))$:
\bel{BVP0}\3n\left\{\2n\ba{ll}
\ns\ds\dot
X(s)=A(s,X(s))-B(s,X(s))R(s,X(s))^{-1}\[B(s,X(s))^TY(s)+S(s,X(s))\],\\
\ns\ds\dot Y(s)=-\[A_x(s,X(s))-\sum_{j=1}^m e_j^T
R(s,X(s))^{-1}\big\{B(s,X(s))^TY(s)+S(s,X(s))\big\}
B^j_x(s,X(s))\]^TY(s)\\
\ns\ds\qq\qq-Q_x(s,X(s))^T+S_x(s,X(s))^TR(s,X(s))^{-1}\[B(s,X(s))^TY(s)+
S(s,X(s))\]\\
\ns\ds\qq\qq-{1\over2}\sum_{j,k=1}^m\[B(s,X(s))^TY(s)+S(s,X(s))\]^T
R(s,X(s))^{-1}e_je_k^TR(s,X(s))^{-1}\\
\ns\ds\qq\qq\qq\qq\cd\[B(s,X(s))^TY(s)+S(s,X(s))\]
R^{jk}_x(s,X(s))^T,\qq s\in[t,T],\\
\ns\ds X(t)=x,\qq Y(T)=G_x(X(T))^T,\ea\right.\ee
where $e_j\in\dbR^m$ is the vector with entry 1 at the $i$-th
position and all other entries are zero. If $(X(\cd),Y(\cd))$ is the
unique solution to the above, then $X(\cd)=X^*(\cd)$ must be the
optimal trajectory and the optimal control $u^*(\cd)$ is given by
(\ref{3.5}).

\rm

\section{Variations of the Cost Functional.}

In the rest of this paper, we let $U=\dbR^m$. In this case,
$\cU[t,T]$ is a Hilbert space whose dual $\cU[t,T]^*$ can be
identified with $\cU[t,T]$ by the Riesz representation theorem. Let
us first make an observation. Define
\bel{F}F(t,x,u(\cd))=D_uJ(t,x;u(\cd)),\qq\forall(t,x,u(\cd))\in[0,T]\times\dbR^n\times\cU[0,T].
\ee
Then $F:[0,T]\times\dbR^n\times\cU[0,T]\to\cU[0,T]^*=\cU[0,T]$. For
any fixed initial pair $(t,x)\in [0,T]\times\dbR^n$, consider the
following equation:
\bel{F(t,x,u)=0}F(t,x,u(\cd))=0.\ee
Under (H1)--(H3), from Proposition 3.1, $u(\cd)\mapsto
J(t,x;u(\cd))$ admits a minimum $u^*(\cd)\equiv
u^*(\cd\,;t,x)\in\cU[t,T]$, i.e.,
\bel{}V(t,x)=J(t,x;u^*(\cd))=\inf_{u(\cd)\in\cU[t,T]}J(t,x;u(\cd)).\ee
Then it is necessary that $u^*(\cd)$ is a solution to equation
(\ref{F(t,x,u)=0}), and
\bel{F_u>0}F_u(t,x;u^*(\cd))=D_{uu}J(t,x;u^*(\cd))\ge0.\ee
Now, suppose $F_u(t,x;u^*(\cd))^{-1}:\cU[t,T]\to\cU[t,T]$ exists and
suppose it is a bounded operator, which, by combining (\ref{F_u>0}),
is equivalent to the following:
\bel{F_u>>0}F_u(t,x;u^*(\cd))=D_{uu}J(t,x;u^*(\cd))\ge\d I,\ee
for some $\d>0$. Then, by implicit function theorem, we have that
$u^*(\cd)\equiv u^*(\cd\,;t,x)$ is differentiable and
\bel{}\ba{ll}
\ns\ds
u^*_{(t,x)}(\cd\,;t,x)=F_u(t,x;u^*(\cd))^{-1}F_{(t,x)}(t,x;u^*(\cd))\\
\ns\ds\qq\qq\q~\equiv-D_{uu}J(t,x;u^*(\cd\,;t,x))^{-1}
[D_uJ]_{(t,x)}(t,x;u^*(\cd\,;t,x)).\ea\ee
Therefore, under (H1)--(H3), as long as
$D_{uu}J(t,x;u^*(\cd\,;t,x))$ is uniformly positive definite,
$(t,x)\mapsto u^*(\cd\,;t,x)$ is differentiable, which implies that
$$V(t,x)\equiv J(t,x;u^*(\cd\,;t,x))$$
is differentiable.

\ms

We now try to find conditions under which (\ref{F_u>>0}) holds. To
this end, let us calculate $D_uJ(t,x;u(\cd))$ and
$D_{uu}J(t,x;u(\cd))$. Denote
$$x\1n=\1n\left(\ba{c}x^1\\ x^2\\ \vdots\\
x^n\ea\right)\2n,~y\1n=\1n\left(\ba{c}y^1\\ y^2\\ \vdots\\
y^n\ea\right)\1n\in\1n\dbR^n,\qq u\1n=\1n\left(\ba{c}u^1\\
u^2\\ \vdots\\ u^m\ea\right)\1n\in\1n\dbR^m,$$
and
$$\ba{ll}
\ns\ds A(t,x)=\left(\ba{c}A^1(t,x)\\ A^2(t,x)\\ \vdots\\
A^n(t,x)\ea\right),\q
B(t,x)\1n=\1n\left(\ba{cccc}B^{11}(t,x)&B^{12}(t,x)&\cds&B^{1m}(t,x)\\
                      B^{21}(t,x)&B^{22}(t,x)&\cds&B^{2m}(t,x)\\
                      \vdots&\vdots&\ddots&\vdots\\
                      B^{n1}(t,x)&B^{n2}(t,x)&\cds&B^{nm}(t,x)\ea\right),\\ [12mm]
\ns\ds S(t,x)=\left(\ba{c}S^1(t,x)\\ S^2(t,x)\\ \vdots\\
S^m(t,x)\ea\right),\q
R(t,x)\1n=\1n\left(\ba{cccc}R^{11}(t,x)&R^{12}(t,x)&\cds&R^{1m}(t,x)\\
                      R^{21}(t,x)&R^{22}(t,x)&\cds&R^{2m}(t,x)\\
                      \vdots&\vdots&\ddots&\vdots\\
                      R^{m1}(t,x)&R^{m2}(t,x)&\cds&R^{mm}(t,x)\ea\right),\\ [10mm]
\ns\ds A^i,B^{ij},S^j,R^{jk}:[0,T]\times\dbR^n\to\dbR,\qq 1\le i\le n,~1\le j,k\le m.\ea$$
Next, we denote
$$\ba{ll}
\ns\ds B^j(t,x)=\left(\ba{c}B^{1j}(t,x)\\ B^{2j}(t,x)\\ \vdots\\
B^{nj}(t,x)\ea\right),~\wt B^i(t,x)=\left(\ba{c}B^{i1}(t,x)\\
B^{i2}(t,x)\\ \vdots\\ B^{im}(t,x)\ea\right),\\ [8mm]
\ns\ds\qq\qq\qq\qq\qq\forall(t,x) \in[0,T]\times\dbR^n,~1\le j\le
m,~1\le i\le n.\ea$$
Then,
$$\left\{\2n\ba{ll}
\ns\ds B(t,x)=(B^1(t,x),B^2(t,x),\cds,B^m(t,x)),\\
\ns\ds B(t,x)^T=(\wt B^1(t,x),\wt B^2(t,x),\cds,\wt
B^n(t,x)),\ea\right.\qq(t,x)\in[0,T]\times\dbR^n.$$

We have the following result.

\ms

\bf Proposition 4.1. \sl Let {\rm(H1)--(H3)} hold. Then for any
$(t,x)\in[0,T)\times\dbR^n$ and $u(\cd)\in\cU[t,T]$,
\bel{DJ}\ba{ll}
\ns\ds
[D_uJ(t,x;u(\cd))](s)=R(s,X(s))u(s)+S(s,X(s))+B(s,X(s))^TY(s),\qq
s\in[t,T],\ea\ee
with $(X(\cd),Y(\cd))$ being the solution to the following decoupled
two-point boundary value problem:
\bel{BVP4.8}\left\{\2n\ba{ll}
\ns\ds\dot X(s)=A(s,X(s))+B(s,X(s))u(s),\\
\ns\ds\dot
Y(s)=-\[A_x(s,X(s))+\sum_{j=1}^mu^j(s)B^j_x(s,X(s))\]^TY(s)
-Q_x(s,X(s))^T\\
\ns\ds\qq\qq-S_x(s,X(s))^Tu(s)-{1\over2}\sum_{j,k=1}^mu^j(s)u^k(s)R^{jk}_x(s,X(s))^T,\\
\ns\ds X(t)=x,\qq Y(T)=G_x(X(T))^T.\ea\right.\ee
Further, for any $v(\cd)\in\cU[t,T]$,
\bel{D^2J}[D_{uu}J(t,x;u(\cd))v(\cd)](s)=R(s,X(s))v(s)+B(s,X(s))^TY_1(s)+\BC(s)X_1(s),\qq
s\in[t,T],\ee
where
\bel{BVP4.10}\left\{\2n\ba{ll}
\ns\ds\left(\ba{c}\dot X_1(s)\\ \dot
Y_1(s)\ea\right)=\left(\ba{cc}\BA(s)&0\\-\BA_1(s)&-\BA(s)^T\ea\right)
\left(\ba{c}X_1(s)\\
Y_1(s)\ea\right)+\left(\ba{c}B(s,X(s))\\-\BC(s)^T\ea\right)v(s),\qq
s\in[t,T],\\
\ns\ds X_1(t)=0,\qq Y_1(T)=G_{xx}(X(T))X_1(T),\ea\right.\ee
with
\bel{4.11}\left\{\2n\ba{ll}
\ns\ds\BA(s)=A_x(s,X(s))+\sum_{j=1}^mu^j(s)B^j_x(s,X(s)),\q\\
\ns\ds\BA_1(s)=\sum_{i=1}^nY^i(s)A^i_{xx}(s,X(s))+Q_{xx}(s,X(s))+\1n{1\over2}\2n\sum_{j,k=1}^mu^j(s)u^k(s)R^{jk}_{xx}(s,X(s))\\
\ns\ds\qq\qq+\2n\sum_{j=1}^mu^j(s)\[\sum_{i=1}^nY^i(s)B^{ij}_{xx}(s,X(s))
\1n+\1n S^j_{xx}(s,X(s))\],\\
\ns\ds\BC(s)=\sum_{j=1}^mu^j(s)R^j_x(s,X(s))+S_x(s,X(s))+\sum_{i=1}^nY^i(s)\wt
B^i_x(s,X(s)).\ea\right.\ee

\ms

\it Proof. \rm Let $(t,x)\in[0,T)\times\dbR^n$ be fixed and
$u(\cd),v(\cd)\in\cU[t,T]$, let
\bel{}\ba{ll}
\ns\ds X(\cd)=X(\cd\,;t,x,u(\cd)),\q X^\e(\cd)=X(\cd\,;t,x,u(\cd)+\e
v(\cd)),\ea\ee
with $\e>0$. Let
\bel{}X_1(\cd)=\lim_{\e\to0}{X^\e(\cd)-X(\cd)\over\e}.\ee
Then
$$\ba{ll}
\ns\ds\dot
X_1(s)=\lim_{\e\to0}\Big\{{A(s,X^\e(s))-A(s,X(s))\over\e}+\sum_{j=1}^mu^j(s)
{B^j(s,X^\e(s))-B^j(s,X(s))\over\e}\Big\}+B(s,X(s))v(s)\\
\ns\ds\qq~=\[A_x(s,X(s))+\sum_{j=1}^mu^j(s)B^j_x(s,X(s))\]X_1(s)+B(s,X(s))v(s).\ea$$
Thus, $X_1(\cd)$ solves the following:
$$\left\{\2n\ba{ll}
\ns\ds\dot X_1(s)=\BA(s)X_1(s)+B(s,X(s))v(s),\qq s\in[t,T],\\
\ns\ds X_1(t)=0,\ea\right.$$
with $\BA(\cd)$ being defined in (\ref{4.11}). We have
\bel{4.8}\ba{ll}
\ns\ds\lan D_uJ(t,x;u(\cd)),v(\cd)\ran=\lim_{\e\to0}{J(t,x;u(\cd)+\e
v(\cd))-J(t,x;u(\cd))\over\e}\\
\ns\ds=\int_t^T\[Q_x(s,X(s))X_1(s)+\lan S(s,X(s)),v(s)\ran+\lan
S_x(s,X(s))X_1(s),u(s)\ran\\
\ns\ds\qq+\lan R(s,X(s))u(s),v(s)\ran+{1\over2}\sum_{j,k=1}^mu^j(s)u^k(s)R^{jk}_x(s,X(s))X_1(s)\]ds+G_x(X(T))X_1(T)\\
\ns\ds=\int_t^T\[\lan
Q_x(s,X(s))^T+S_x(s,X(s))^Tu(s)+{1\over2}\sum_{j,k=1}^mu^j(s)u^k(s)R^{jk}_x(s,X(s))^T,X_1(s)\ran\\
\ns\ds\qq+\lan
S(s,X(s))+R(s,X(s))u(s),v(s)\ran\]ds+G_x(X(T))X_1(T).\ea\ee
Let $(X(\cd),Y(\cd))$ be the solution to (\ref{BVP4.8}). Then (note
(\ref{4.11}))
$$\ba{ll}
\ns\ds{d\over ds}\lan Y(s),X_1(s)\ran=\lan \dot Y(s),X_1(s)\ran+\lan
Y(s),\BA(s)X_1(s)\ran+\lan Y(s),B(s,X(s))v(s)\ran\\
\ns\ds=\lan\dot Y(s)+\BA(s)^TY(s),X_1(s)\ran+\lan
B(s,X(s))^TY(s),v(s)\ran.\ea$$
Noting $X_1(t)=0$, one has
$$\ba{ll}
\ns\ds G_x(X(T))X_1(T)\3n\2n&\ds=\lan Y(T),X_1(T)\ran\\
\ns&\ds=\int_t^T\Big\{\lan\dot Y(s)+\BA(s)^TY(s),X_1(s)\ran+\lan
B(s,X(s))^TY(s),v(s)\ran\Big\}ds.\ea$$
Consequently,
$$\ba{ll}
\ns\ds\lan D_uJ(t,x;u(\cd)),v(\cd)\ran=\int_t^T\Big\{\lan\dot
Y(s)+\BA(s)^TY(s)\\
\ns\ds\qq+Q_x(s,X(s))^T+S_x(s,X(s))^Tu(s)+{1\over2}\sum_{j,k=1}^mu^j(s)u^k(s)R^{jk}_x(s,X(s))^T,X_1(s)\ran\\
\ns\ds\qq+\lan
B(s,X(s))^TY(s)+S(s,X(s))+R(s,X(s))u(s),v(s)\ran\Big\}ds\\
\ns\ds=\int_t^T\lan
R(s,X(s))u(s)+S(s,X(s))+B(s,X(s))^TY(s),v(s)\ran ds.\ea$$
This proves (\ref{DJ}).

\ms

Next, we calculate $D_{uu}J(t,x;u(\cd))$. To this end, for any
$\e\in(0,1)$, let $(X^\e(\cd),Y^\e(\cd))$ be the solution to the
following:
\bel{}\left\{\2n\ba{ll}
\ns\ds\dot X^\e(s)=A(s,X^\e(s))+B(s,X^\e(s))\big[u(s)+\e v(s)],\\
\ns\ds\dot Y^\e(s)=-\[A_x(s,X^\e(s))+\sum_{j=1}^m[u^j(s)+\e v^j(s)]
B^j_x(s,X^\e(s))\]^TY^\e(s)\\
\ns\ds\qq\qq-Q_x(s,X^\e(s))^T-S_x(s,X^\e(s))^T[u(s)+\e v(s)]\\
\ns\ds\qq\qq-{1\over2}\sum_{j,k=1}^m[u^j(s)+\e v^j(s)][u^k(s)+\e
v^k(s)]R^{jk}_x(s,X^\e(s))^T,\\
\ns\ds X^\e(t)=x,\qq Y^\e(T)=G_x(X^\e(T))^T.\ea\right.\ee
Then
$$\ba{ll}
\ns\ds[D_uJ(t,x;u(\cd)\1n+\1n\e v(\cd))](s)\1n=\1n
R(s,X^\e(s))[u(s)\1n+\1n\e
v(s)]\1n+\1n S(s,X^\e(s))\1n+\1n B(s,X^\e(s))^TY^\e(s),\\
\ns\ds\qq\qq\qq\qq\qq\qq\qq\qq\qq\qq\qq\qq\qq\qq s\in[t,T].\ea$$
Hence,
$$\ba{ll}
\ns\ds[D_{uu}J(t,x;u(\cd))v(\cd)](s)=\lim_{\e\to0}{[D_uJ(t,x;u(\cd)+\e
v(\cd))](s)-[D_uJ(t,x;u(\cd))](s)\over\e}\\
\ns\ds=\lim_{\e\to0}\Big\{R(s,X^\e(s))v(s)+{R(s,X^\e(s))-R(s,X(s))\over\e}u(s)
+{S(x,X^\e(s))-S(s,X(s))\over\e}\\
\ns\ds\qq\qq+B(s,X^\e(s))^T{Y^\e(s)-Y(s)\over\e}+{B(s,X^\e(s))^T-B(s,X(s))^T
\over\e}Y(s)\Big\}\\
\ns\ds=R(s,X(s))v(s)+\sum_{j=1}^mu^j(s)R^j_x(s,X(s))X_1(s)+S_x(s,X(s))X_1(s)\\
\ns\ds\qq+B(s,X(s))^TY_1(s)+\sum_{i=1}^nY^i(s)\wt
B^i_x(s,X(s))X_1(s)\\
\ns\ds=R(s,X(s))v(s)+B(s,X(s))^TY_1(s)\\
\ns\ds\qq\qq+\[\sum_{j=1}^mu^j(s)R^j_x(s,X(s))+S_x(s,X(s))+\sum_{i=1}^nY^i(s)\wt
B^i_x(s,X(s))\]X_1(s)\\
\ns\ds\equiv R(s,X(s))v(s)+B(s,X(s))^TY_1(s)+\BC(s)X_1(s),\ea$$
where
\bel{}Y_1(s)=\lim_{\e\to0}{Y^\e(s)-Y(s)\over\e}\,,\ee
and $\BC(\cd)$ is defined in (\ref{4.11}). Then to complete the
proof, we need only to derive the equation for $Y_1(\cd)$. First of
all,
\bel{}\ba{ll}
\ns\ds Y_1(T)=\lim_{\e\to0}{Y^\e(T)-Y(T)\over\e}=\lim_{\e\to0}{G_x(X^\e(T))^T
-G_x(X(T))^T\over\e}\\
\ns\ds\qq\q=G_{xx}(X(T))\left[\lim_{\e\to0}{X^\e(T)-X(T)\over\e}\right]=G_{xx}(X(T))X_1(T).\\
\ea\ee
Next,
$$\ba{ll}
\ns\ds\dot
Y_1(s)=\lim_{\e\to0}{\dot Y^\e(s)-\dot Y(s)\over\e}\\
\ns\ds\qq=-\lim_{\e\to0}\Big\{A_x(s,X^\e(s))^T
{Y^\e(s)-Y(s)\over\e}+\sum_{i=1}^nY^i{A^i_x(s,X^\e(s))^T-A^i_x(s,X(s))^T\over\e}\\
\ns\ds\qq\qq+\sum_{j=1}^mv^j(s)B^j_x(s,X^\e(s))^TY^\e(s)
+\sum_{j=1}^mu^j(s)B^j_x(s,X^\e(s))^T
{Y^\e(s)-Y(s)\over\e}\\
\ns\ds\qq\qq+\sum_{j=1}^m\sum_{i=1}^nu^j(s){B^{ij}_x(s,X^\e(s))^T
-B^{ij}_x(s,X(s))^T\over\e}Y^i(s)+{Q_x(X^\e(s))^T-Q_x(X(s))^T\over\e}\\
\ns\ds\qq\qq+S_x(s,X^\e(s))^Tv(s)+{S_x(s,X^\e(s))^T-S_x(s,X(s))^T\over\e}u(s)\\
\ns\ds\qq\qq+{1\over2}\2n\sum_{j,k=1}^m\2n
u^j(s)u^k(s){R^{jk}_x(s,X^\e(s))^T\2n-\1n
R^{jk}_x(s,X(s))^T\over\e}\1n
+\2n\sum_{j,k=1}^mu^j(s)v^k(s)R^{jk}_x(s,X^\e(s))^T \Big\}.\ea$$
Hence,
$$\ba{ll}
\ns\ds\dot
Y_1(s)=-A_x(s,X(s))^TY_1(s)-\sum_{i=1}^nY^i(s)A^i_{xx}(s,X(s))X_1(s)
-\sum_{j=1}^mv^j(s)B^j_x(s,X(s))^TY(s)\\
\ns\ds\qq\qq\q-\sum_{j=1}^mu^j(s) B^j_x(s,X(s))^T
Y_1(s)-\sum_{j=1}^m\sum_{i=1}^nu^j(s)Y^i(s)B^{ij}_{xx}(s,X(s))X_1(s)\\
\ns\ds\qq\qq\q-Q_{xx}(s,X(s))X_1(s)-S_x(s,X(s))^Tv(s)-\sum_{j=1}^mu^j(s)S^j_{xx}(s,X(s))X_1(s)\\
\ns\ds\qq\qq-{1\over2}\sum_{j,k=1}^mu^j(s)u^k(s)R^{jk}_{xx}(s,X(s))X_1(s)-\sum_{j,k=1}^mu^j(s)v^k(s)
R^{jk}_x(s,X(s))^T\\
\ns\ds\qq=-\[A_x(s,X(s))+\sum_{j=1}^mu^j(s)B^j_x(s,X(s))\]^T
Y_1(s)\\
\ns\ds\qq\qq-\[\sum_{i=1}^nY^i(s)A^i_{xx}(s,X(s))
+\sum_{j=1}^m\sum_{i=1}^nu^j(s)Y^i(s)B^{ij}_{xx}(s,X(s))
+Q_{xx}(s,X(s))\\
\ns\ds\qq\qq\qq+\sum_{j=1}^mu^j(s)S^j_{xx}(s,X(s))
+{1\over2}\sum_{j,k=1}^mu^j(s)u^k(s)R^{jk}_{xx}(s,X(s))\]
X_1(s)\\
\ns\ds\qq\qq-\[S_x(s,X(s))+\sum_{i=1}^nY^i(s)\wt
B^i_x(s,X(s))+\sum_{j=1}^mu^j(s)R^j_x(s,X(s))\]^Tv(s)\\
\ns\ds\qq=-\BA(s)^TY_1(s)-\BA_1(s)X_1(s)-\BC(s)^Tv(s),\ea$$
where $\BA(\cd)$, $\BA_1(\cd)$, and $\BC(\cd)$ are given by
(\ref{4.11}). Thus, $(X_1(\cd),Y_1(\cd))$ solves (\ref{BVP4.10}).
\endpf

\ms

Note that for given $(t,x)\in[0,T)\times\dbR^n$ and
$u(\cd)\in\cU[t,T]$, both $X_1(\cd)$ and $Y_1(\cd)$ depend on
$v(\cd)$. It will be desirable to have a representation of
$[D_{uu}J(t,x;u(\cd))v(\cd)]$ explicitly in terms of $v(\cd)$. The
following is such a result.

\ms

\bf Proposition 4.2. \sl Let {\rm(H1)--(H3)} hold. For any
$(t,x)\in[0,T)\times\dbR^n$ and $u(\cd)\in\cU[t,T]$, let $\BA(\cd)$,
$\BA_1(\cd)$, and $\BC(\cd)$ be defined by $(\ref{4.11})$. Then
\bel{D^2J(u)v}[D_{uu}J(t,x;u(\cd))v(\cd)](s)=R(s,X(s))v(s)
+\int_t^T\BF(s,r)v(r)ds,\q\forall v(\cd)\in\cU[t,T],\ee
where
\bel{BF}\ba{ll}
\ns\ds\BF(s,r)=B(s,X(s))^T\1n\F_\BA(T,s)^TG_{xx}(X(T))
\F_\BA(T,r)B(r,X(r))\\
\ns\ds\qq\qq+\int_{s\vee
r}^TB(s,X(s))^T\F_\BA(r',s)^T\2n\BA_1(r')\F_\BA(r',r)B(r,X(r))dr'\\
\ns\ds\qq\qq+\BC(s)\F_\BA(s,r)B(r,X(r))I_{[t,s]}(r)+B(s,X(s))^T\F_\BA(r,s)^T\BC(r)^T
I_{[s,T]}(r),\ea\ee
and $\F_\BA(\cd\,,\cd)$ is the fundamental matrix of $\BA(\cd)$,
i.e., for any $\t\in[t,T)$,
$$\left\{\2n\ba{ll}
\ns\ds{d\over ds}\F_\BA(s,\t)=\BA(s)\F_\BA(s,\t),\qq s\in[\t,T],\\
\ns\ds\F_\BA(\t,\t)=I.\ea\right.$$

\it Proof. \rm Let $\F_\BA(\cd\,,\cd)$ be the fundamental matrix of
$\BA(\cd)$. Then
$$X_1(s)=\int_t^s\F_\BA(s,r)B(r,X(r))v(r)dr,\qq s\in[t,T],$$
and
$$\ba{ll}
\ns\ds Y_1(s)=\F_\BA(T,s)^TG_{xx}(X(T))X_1(T)+\int_s^T\F_\BA(r,s)^T
\[\BA_1(r)X_1(r)
+\BC(r)^Tv(r)\]dr\\
\ns\ds\qq=\F_\BA(T,s)^TG_{xx}(X(T))\int_t^T\F_\BA(T,r)B(r,X(r))v(r)dr\\
\ns\ds\qq\qq+\1n\int_s^T\1n\F_\BA(r,s)^T\1n\BA_1(r)\1n\int_t^r\2n
\F_\BA(r,r')B(r',X(r'))v(r')dr'dr
\1n+\2n\int_s^T\2n\F_\BA(r,s)^T\1n\BC(r)^T\1n v(r)dr\\
\ns\ds\qq=\F_\BA(T,s)^TG_{xx}(X(T))\int_t^T\F_\BA(T,r)B(r,X(r))v(r)dr\\
\ns\ds\qq\qq+\2n\int_t^T\2n\[\int_{s\vee
r}^T\2n\F_\BA(r',s)^T\1n\BA_1(r')\F_\BA(r',r)dr'\]B(r,X(r))v(r)dr
\1n+\2n\int_s^T\2n\F_\BA(r,s)^T\1n\BC(r)^T\1n v(r)dr\ea$$
Hence,
$$\ba{ll}
\ns\ds[D_{uu}J(t,x;u(\cd))v(\cd)](s)=R(s,X(s))v(s)+B(s,X(s))^TY_1(s)+\BC(s)X_1(s)\\
\ns\ds=R(s,X(s))v(s)+\BC(s)\int_t^s\F_\BA(s,r)B(r,X(r))v(r)dr\\
\ns\ds\q+B(s,X(s))^T\[\F_\BA(T,s)^TG_{xx}(X(T))\int_t^T\F_\BA(T,r)B(r,X(r))v(r)dr\\
\ns\ds\qq\qq+\2n\int_t^T\2n\(\int_{s\vee
r}^T\2n\F_\BA(r',s)^T\1n\BA_1(r')\F_\BA(r',r)dr'\)B(r,X(r))v(r)dr
\1n+\2n\int_s^T\2n\F_\BA(r,s)^T\1n\BC(r)^T\1n v(r)dr\]\\
\ns\ds=R(s,X(s))v(s)+\int_t^T\BF(s,r)v(r)ds,\ea$$
proving (\ref{D^2J(u)v}). \endpf

\ms

We note that $\BF(s,r)$ is depending on the given $u(\cd)$ and is
independent of $v(\cd)$.

\section{Invertibility of $D_{uu}J(t,x;u(\cd))$.}

Having calculated $D_{uu}J(t,x;u(\cd))$, we now would like to look
at conditions under which it admits a bounded inverse. The following
is a general result whose proof is straightforward.

\ms

\bf Proposition 5.1. \sl Let {\rm(H1)--(H3)} hold and let
$(t,x)\in[0,T)\times\dbR^n$, $u(\cd)\in\cU[t,T]$ be given. Define
$\BF(\cd\,,\cd)$ by $(\ref{BF})$, with $(X(\cd),Y(\cd))$ being the
solution to $(\ref{BVP4.8})$. Then $D_{uu}J(t,x;u(\cd))$ admits a
bounded inverse operator if and only if for any $w(\cd)\in\cU[t,T]$,
the following second kind Fredholm integral equation is well-posed:
\bel{Fredholm}w(s)=R(s,X(s))v(s)+\int_t^T\BF(s,r)v(r)dr,\qq
s\in[t,T].\ee
A sufficient condition for the above is
\bel{<1}|R(s,X(s))^{-1}\BF(s,r)|\le\a<{1\over T}\,,\qq
s,r\in[0,T].\ee

\rm

Practically, to use the above result, we need to first solve a
(decoupled) two-point boundary value problem (\ref{BVP4.8}) to get
$(X(\cd),Y(\cd))$. Then calculate $\BA(\cd)$, $\BA_1(\cd)$ and
$\BC(\cd)$, etc., followed by $\F_\BA(\cd\,,\cd)$. Next, construct
$\BF(\cd\,,\cd)$ and then check see if the Fredholm integral
equation (\ref{Fredholm}) is well-posed or sufficiently look at if
(\ref{<1}) holds. Apparently, some more direct sufficient conditions
are desirable for $D_{uu}J(t,x;u(\cd))$ to be invertible. Recall
from the previous section that the invertibility of
$D_{uu}J(t,x;u(\cd))$ is equivalent to the uniform positive
definiteness (see (\ref{F_u>>0})):
\bel{5.3}D_{uu}J(t,x;u(\cd))\ge\d I,\ee
for some $\d>0$. Thus, we now would like to look for some sufficient
conditions under which (\ref{5.3}) is satisfied. To approach this,
we first present the following proposition.

\ms

\bf Proposition 5.2. \sl Let {\rm(H1)--(H3)} hold. Let
$(t,x)\in[0,T)\times\dbR^n$, and $u(\cd)\in\cU[t,T]$ be given. Let
$(X(\cd),Y(\cd))$ be the solution to $(\ref{BVP4.8})$ and
$\BA(\cd)$, $\BA_1(\cd)$, and $\BC(\cd)$ be defined by
$(\ref{4.11})$. Then
\bel{5.4}\ba{ll}
\ns\ds\int_t^T\lan\,[D_{uu}J(t,x;u(\cd))v(\cd)](s),v(s)\ran ds
=\int_t^T\3n\lan R(s,X(s))v(s),v(s)\ran ds\\
\ns\ds\q+\lan G_{xx}(X(T))\int_t^T\F_\BA(T,r)B(r,X(r))v(r)dr,
\int_t^T\F_\BA(T,r)B(r,X(r))v(r)dr\ran\\
\ns\ds\q+\int_t^T\lan\BA_1(s)\int_t^s\F_\BA(s,r)B(r,X(r))v(r)dr,
\int_t^s\F_\BA(s,r)B(r,X(r))v(r)dr\ran ds\\
\ns\ds\q+2\int_t^T\lan\BC(s)\int_t^s\F_\BA(s,r)B(r,X(r))v(r)dr,v(s)\ran
ds.\ea\ee
Further, suppose $\bar G\in\dbS^n_+$ and $\bar Q:[0,T]\to\dbS^n_+$
such that for some $\a\in(0,1)$,
\bel{G+G>0}\left\{\ba{ll}
\ns\ds G_{xx}(X(T))+\bar G\ge0,\\
\ns\ds\BA_1(s)+\bar Q(s)-\a^{-1}\BC(s)^TR(s,X(s))^{-1}\BC(s)\ge0,\qq
s\in[t,T],\ea\right.\ee
and
\bel{}(1-\a)R(s,X(s))-\big[\h G(t)+\h Q(s,t)\big]I\ge\d I,\qq
s\in[t,T],\ee
for some $\d>0$, with
\bel{5.6}\left\{\2n\ba{ll}
\ns\ds\h G(t)=\[\int_t^T\int_t^T|B(s,X(s))^T\F_\BA(T,s)^T\bar
G\F_\BA(T,r)B(r,X(r))|^2drds\]^{1\over2},\\
\ns\ds\h
Q(s,t)=\int_s^T\[\int_t^\t\int_t^\t|B(r,X(r))^T\F_\BA(\t,r)^T\bar
Q(\t)\F_\BA(\t,r')B(r',X(r'))|^2dr'dr\]^{1\over2}d\t,\ea\right.\ee
then
\bel{D2J>d}D_{uu}J(t,x;u(\cd))\ge\d I.\ee

\ms

\it Proof. \rm Let $(t,x)\in[0,T)\times\dbR^n$ and
$u(\cd)\in\cU[t,T]$ be given. We have
$$\ba{ll}
\ns\ds\int_t^T\int_t^T\lan\BF(s,r)v(r),v(s)\ran
drds\\
\ns\ds=\lan
G_{xx}(X(T))\int_t^T\F_\BA(T,r)B(r,X(r))v(r)dr,\int_t^T\F_\BA(T,r)B(r,X(r))v(r)dr\ran\\
\ns\ds\q+\int_t^T\lan\BA_1(s)\int_t^s\F_\BA(s,r)B(r,X(r))v(r)dr,
\int_t^s\F_\BA(s,r)B(r,X(r))v(r)dr\ran ds\\
\ns\ds\q+\int_t^T\3n\int_t^s\lan\BC(s)\F_\BA(s,r)B(r,X(r))v(r),v(s)\ran
drds\\
\ns\ds\q+\int_t^T\3n\int_s^T\lan
B(s,X(s))^T\F_\BA(r,s)^T\BC(r)^Tv(r),v(s)\ran drds\ea$$
$$\ba{ll}
\ns\ds=\lan
G_{xx}(X(T))\int_t^T\F_\BA(T,r)B(r,X(r))v(r)dr,\int_t^T\F_\BA(T,r)B(r,X(r))v(r)dr\ran\\
\ns\ds\q+\int_t^T\lan\BA_1(s)\int_t^s\F_\BA(s,r)B(r,X(r))v(r)dr,
\int_t^s\F_\BA(s,r)B(r,X(r))v(r)dr\ran ds\\
\ns\ds\q+2\int_t^T\lan\BC(s)\int_t^s\F_\BA(s,r)B(r,X(r))v(r)dr,v(s)\ran
ds.\ea$$
This proves (\ref{5.4}). From this, one further has
$$\ba{ll}
\ns\ds\int_t^T\lan[D_{uu}J(t,x;u(\cd))v(\cd)](s),v(s)\ran
ds=\int_t^T\lan R(s,X(s))v(s),v(s)\ran ds\\
\ns\ds\q+\lan
G_{xx}(X(T))\int_t^T\F_\BA(T,r)B(r,X(r))v(r)dr,\int_t^T\F_\BA(T,r)B(r,X(r))v(r)dr\ran\\
\ns\ds\q+\int_t^T\lan\BA_1(s)\int_t^s\F_\BA(s,r)B(r,X(r))v(r)dr,
\int_t^s\F_\BA(s,r)B(r,X(r))v(r)dr\ran ds\\
\ns\ds\q+2\int_t^T\lan\BC(s)\int_t^s\F_\BA(s,r)B(r,X(r))v(r)dr,v(s)\ran
ds\\
\ns\ds=\int_t^T\lan(1-\a)R(s,X(s))v(s),v(s)\ran
ds+\int_t^T\(\big|\a^{1\over2}R(s,X(s))^{1\over2}v(s)\big|^2ds\\
\ns\ds\q+2\lan\a^{-{1\over2}}R(s,X(s))^{-{1\over2}}\BC(s)\int_t^s\F_\BA(s,r)
B(r,X(r))v(r)dr,\a^{1\over2}R(s,X(s))^{1\over2}v(s)\ran\\
\ns\ds\q+\big|\a^{-{1\over2}}R(s,X(s))^{-{1\over2}}\BC(s)\int_t^s\F_\BA(s,r)
B(r,X(r))v(r)dr\big|^2\)ds\\
\ns\ds\q+\lan G_{xx}(X(T))\int_t^T\F_\BA(T,r)B(r,X(r))v(r)dr,
\int_t^T\F_\BA(T,r)B(r,X(r))v(r)dr\ran\\
\ns\ds\q+\int_t^T\lan
\[\BA_1(s)-\a^{-1}\BC(s)^TR(s,X(s))^{-1}\BC(s)\]\int_t^s\F_A(s,r)B(r,X(r))v(r)dr,\\
\ns\ds\qq\qq\qq\qq\qq\qq\qq\qq\qq\qq\qq\int_t^s\F_A(s,r)B(r,X(r))v(r)dr\ran ds\\
\ns\ds=\int_t^T\lan(1-\a)R(s,X(s))v(s),v(s)\ran ds\\
\ns\ds\q+\int_t^T\Big|\,\a^{1\over2}R(s,X(s))^{1\over2}v(s)
+\a^{-{1\over2}}R(s,X(s))^{-{1\over2}}\BC(s)\int_t^s\F_\BA(s,r)B(r,X(r))v(r)dr\Big|^2ds\\
\ns\ds\q+\lan\big[G_{xx}(X(T))+\bar
G\big]\int_t^T\F_\BA(T,r)B(r,X(r))v(r)dr,
\int_t^T\F_\BA(T,r)B(r,X(r))v(r)dr\ran\\
\ns\ds\q+\2n\int_t^T\3n\lan\1n\big[\BA_1(s)\1n+\1n\bar
Q(s)\1n-\1n\a^{-1}\BC(s)^T\1n R(s,X(s))^{-1}\BC(s)\1n\big]\2n
\int_t^s\3n\F_A(s,r)B(r,X(r))v(r)dr,\\
\ns\ds\qq\qq\qq\qq\qq\qq\qq\qq\qq\qq\qq\int_t^s\3n\F_A(s,r)B(r,X(r))v(r)dr\1n\ran\1n ds\\
\ns\ds\q-\lan\bar
G\int_t^T\F_\BA(T,r)B(r,X(r))v(r)dr,\int_t^T\F_\BA(T,r)B(r,X(r))v(r)dr\ran\\
\ns\ds\q-\2n\int_t^T\3n\lan\bar Q(s)\2n
\int_t^s\3n\F_\BA(s,r)B(r)v(r)dr,
\int_t^s\3n\F_\BA(s,r)B(r)v(r)dr\1n\ran\1n ds\\
\ns\ds\ge\1n\int_t^T\2n\lan(1\1n-\1n\a)R(s,X(s))v(s),v(s)\ran
ds\1n-\1n\lan\bar
G\2n\int_t^T\2n\F_\BA(T,r)B(r,X(r))v(r)dr,\int_t^T\2n\F_\BA(T,r)B(r,X(r))v(r)dr\ran\\
\ns\ds\q-\2n\int_t^T\3n\lan\bar Q(s)\2n
\int_t^s\3n\F_\BA(s,r)B(r)v(r)dr,
\int_t^s\3n\F_\BA(s,r)B(r)v(r)dr\1n\ran\1n ds.\ea$$
Note that
$$\ba{ll}
\ns\ds\lan\bar
G\int_t^T\F_\BA(T,r)B(r,X(r))v(r)dr,\int_t^T\F_\BA(T,r)B(r,X(r))v(r)dr\ran\\
\ns\ds=\int_t^T\lan\int_t^TB(s,X(s))^T\F_\BA(T,s)^T\bar
G\F_\BA(T,r)B(r,X(r))v(r)dr,v(s)\ran ds\\
\ns\ds\le\(\int_t^T\Big|\int_t^TB(s,X(s))^T\F_\BA(T,s)^T\bar
G\F_\BA(T,r)B(r,X(r))v(r)dr\Big|^2ds\)^{1\over2}\(\int_t^T|v(s)|^2ds\)^{1\over2}\\
\ns\ds\le\[\int_t^T\(\int_t^T|B(s,X(s))^T\F_\BA(T,s)^T\bar
G\F_\BA(T,r)B(r,X(r))|^2drds\)\(\int_t^T|v(r)|^2dr\)\]^{1\over2}
\(\int_t^T|v(r)|^2dr\)^{1\over2}\\
\ns\ds=\[\int_t^T\2n\int_t^T|B(s,X(s))^T\F_\BA(T,s)^T\bar
G\F_\BA(T,r)B(r,X(r))|^2drds\]^{1\over2}\2n\int_t^T\2n|v(r)|^2dr\equiv\2n
\int_t^T\2n\lan\h G(t)v(s),v(s)\ran ds,\ea$$
and similarly,
$$\ba{ll}
\ns\ds\int_t^T\lan\bar Q(s)\int_t^s\F_\BA(s,r)B(r,X(r)))v(r)dr,
\int_t^s\F_\BA(s,r)B(r,X(r))v(r)dr\ran ds\\
\ns\ds=\int_t^T\int_t^s\lan\int_t^sB(r',X(r'))^T\F_\BA(s,r')\bar
Q(s)\F_\BA(s,r)B(r,X(r))v(r)dr,v(r')\ran dr'
ds\\
\ns\ds\le\int_t^T\[\int_t^s\int_t^s|B(r,X(r))^T\F_\BA(s,r)^T\bar
Q(s)\F_A(s,r')B(r',X(r'))|^2dr'dr\]^{1\over2}\int_t^s|v(r)|^2drds\\
\ns\ds=\int_t^T\Big\{\int_\t^T\[\int_t^s\int_t^s|B(r,X(r))^T\F_\BA(t,x)^T\bar
Q(s)\F_\BA(s,r')B(r',X(r')|^2dr'dr\]^{1\over2}ds\Big\}|v(\t)|^2d\t\\
\ns\ds\equiv\int_t^T\h Q(s,t)|v(s)|^2ds,\ea$$
where $\h G(t)$ and $\h Q(s,t)$ are given by (\ref{5.6}).
Consequently,
$$\ba{ll}
\ns\ds\int_t^T\lan\,[D_{uu}J(t,x;u(\cd))v(\cd)](s),v(s)\ran ds\\
\ns\ds\ge\1n\int_t^T\2n\lan(1\1n-\1n\a)R(s,X(s))v(s),v(s)\ran
ds\1n-\1n\lan\bar
G\2n\int_t^T\2n\F_\BA(T,r)B(r,X(r))v(r)dr,\int_t^T\2n\F_\BA(T,r)B(r,X(r))v(r)dr\ran\\
\ns\ds\q-\2n\int_t^T\3n\lan\bar Q(s)\2n
\int_t^s\3n\F_\BA(s,r)B(r,X(r))v(r)dr,
\int_t^s\3n\F_\BA(s,r)B(r,X(r))v(r)dr\1n\ran\1n ds\\
\ns\ds\ge\int_t^T\lan\((1-\a)R(s,X(s))-[\h G(t)+\h
Q(s,t)]I\)v(s),v(s)\ran ds.\ea$$
Hence, (\ref{D2J>d}) follows. \endpf

\ms

Let us point out that we actually do not need the invertibility of
$D_{uu}J(t,x;u(\cd))$ for all $u(\cd)\in\cU[t,T]$, instead, it will
be enough for us to have the invertibility of $D_{uu}J(t,x;u(\cd))$
for $u(\cd)$ given by (see (\ref{3.5}))
$$u(s)=-R(s,X(s))^{-1}\[B(s,X(s))^TY(s)+S(s,X(s))\],\qq s\in[t,T],$$
with $(X(\cd),Y(\cd))$ being any solution to (\ref{BVP0}) whose
existence is guaranteed by the existence of optimal controls and the
Pontryagin's minimum principle.

\ms

We now look at some interesting cases.

\ms

\subsection{Linear quadratic case.}

\rm Let
\bel{LQ}\left\{\2n\ba{ll}
\ns\ds A(t,x)=A(t)x,\q B(t,x)=B(t),\q Q(t,x)={1\over2}\lan
Q(t)x,x\ran,\\
\ns\ds S(t,s)=S(t)x,\q R(t,x)=R(t),\q G(x)={1\over2}\lan
Gx,x\ran.\ea\right.\ee
This is a classical LQ case. In this case,
$$\BA(s)=A(s),\q\BA_1(s)=Q(s),\q\BC(s)=S(s),\qq s\in[0,T].$$
Then (\ref{G+G>0}) holds if
$$\left\{\ba{ll}
\ns\ds G+\bar G\ge0,\\
\ns\ds Q(s)+\bar Q(s)-\a^{-1}S(s)^TR(s)^{-1}S(s)\ge0,\qq
s\in[0,T].\ea\right.$$
for some $\bar G\in\dbS_+^n$, $\bar Q:[0,T]\to\dbS_+^n$, and
$\a\in(0,1)$. In this case, $\F_A(\cd\,,\cd)$, the fundamental
matrix of $A(\cd)$, is independent of $u(\cd)$, $X(\cd)$ and
$Y(\cd)$. Consequently,
\bel{5.9}\left\{\2n\ba{ll}
\ns\ds\h G(t)=\[\int_t^T\int_t^T|B(s)^T\F_A(T,s)^T\bar
G\F_A(T,r)B(r)|^2drds\]^{1\over2},\\
\ns\ds\h Q(s,t)=\int_s^T\[\int_t^\t\int_t^\t|B(r)^T\F_A(\t,r)^T\bar
Q(\t)\F_A(\t,r')B(r')|^2dr'dr\]^{1\over2}d\t.\ea\right.\ee
are independent of $u(\cd)$, $X(\cd)$, and $Y(\cd)$. Then, by
Proposition 5.2, we obtain the following result.

\ms

\bf Proposition 5.3. \sl Let {\rm (H1)--(H3)} and $(\ref{LQ})$ hold.
Suppose there exist $\a\in(0,1)$, and $\bar G\in\dbS_+^n$, $\bar
Q:[t,T]\to\dbS_+^n$ such that
\bel{}Q(s)+\bar Q(s)-\a^{-1}S(s)^TR(s)^{-1}S(s)\ge0,\q G+\bar
G\ge0,\q s\in[t,T],\ee
and
\bel{}(1-\a)R(s)-[\h G(t)+\h Q(s,t)]I\ge\d I,\qq s\in[t,T],\ee
for some $\d>0$, where $\h G(t)$ and $\h Q(s,t)$ are defined by
$(\ref{5.9})$. Then
$$D_{uu}J(t,x;u(\cd))\ge\d I.$$

\ms

\rm

We point out that under the following classical conditions for LQ
problems:
\bel{classical}R(s)\ge\d I,\q Q(s)-S(s)^TR(s)^{-1}S(s)\ge0,\q
G\ge0,\ee
we need only take
$$\a=0,\q\bar G=0,\q\bar Q(\cd)=0.$$
Therefore, the above result covers the classical LQ problem.
Further, Proposition 5.3 shows that for LQ problems, the failure of
the last two conditions in (\ref{classical}) can be compensated by
the sufficient positive definiteness of $R(s)$. On the other hand,
we see that due to the nature of LQ problem, the positive
definiteness of $D_{uu}J(t,x;u(\cd))$ obtained above is
automatically uniform in $u(\cd)$.

\ms

\subsection{Linear semi-convex case.}

Let us first assume the following:
\bel{convex}\left\{\2n\ba{ll}
\ns\ds A(t,x)=A(t)x,\q B(t,x)=B(t),\q S(t,x)=0,\q
R(t,x)=R(t)\ge\d I,\\
\ns\ds x\mapsto Q(t,x),\q x\mapsto G(x)\q\hb{are
convex}.\ea\right.\ee
In the above case, we have a linear state equation and a convex cost
functional. This is a natural generalization of LQ case and we refer
to it as {\it linear-convex} problem. Such kind of problems were
carefully studied in \cite{You 1987, You 1997} by means of the
so-called quasi-Riccati equation.

\ms

Note that under (\ref{convex}), it is straightforward that
$u(\cd)\mapsto J(t,x;u(\cd))$ is uniformly convex. In our framework,
one has
$$\BA(s)=A(s),\q\BA_1(s)=Q_{xx}(s,X(s)),\q\BC(s)=0,\qq
s\in[0,T].$$
Then
$$\ba{ll}
\ns\ds\int_t^T\lan[D_{uu}J(t,x;u(\cd))v(\cd)](s),v(s)\ran ds\\
\ns\ds=\int_t^T\3n\lan R(s)v(s),v(s)\ran ds+\lan
G_{xx}(X(T))\int_t^T\F_A(T,r)B(r)v(r)dr,
\int_t^T\F_A(T,r)B(r)v(r)dr\ran\\
\ns\ds\q+\int_t^T\lan Q_{xx}(s,X(s))\int_t^s\F_A(s,r)B(r)v(r)dr,
\int_t^s\F_A(s,r)B(r)v(r)dr\ran ds\\
\ns\ds\ge\int_t^T\3n\lan R(s)v(s),v(s)\ran
ds\ge\d\int_t^T|v(s)|^2ds,\ea$$
proving the uniform convexity of the map $u(\cd)\mapsto
J(t,x;u(\cd))$.

\ms

We can actually do a little bit more. Here is the result.

\ms

\bf Proposition 5.4. \sl Let {\rm(H1)--(H3)} hold such that for some
$\bar Q:[0,T]\to\dbS_+^n$ and $\bar G\in\dbS_+^n$,
\bel{5.14}A(t,x)=A(t)x,\q B(t,x)=B(t),\q S(t,x)=S(t)x,\q
R(t,x)=R(t),\ee
and
\bel{5.15}\left\{\ba{ll}
\ns\ds Q_{xx}(s,x)+\bar
Q(s)-\a^{-1}S(s)^TR(s)^{-1}S(s)\ge0,\qq(s,x)\in[t,T]\times\dbR^n,\\
\ns\ds G_{xx}(x)+\bar G\ge0,\q s\in[t,T],\ea\right.\ee
and
\bel{5.16}(1-\a)R(s)-[\h G(t)+\h Q(s,t)]I\ge\d I,\qq s\in[t,T],\ee
for some $\d>0$, where $\h G(t)$ and $\h Q(s,t)$ are defined by
$(\ref{5.9})$. Then $D_{uu}J(t,x;u(\cd))$ is uniformly positive
definite.

\ms

\it Proof. \rm We note that under our conditions, one has
$$\BA(s)=A(s),\q\BA_1(s)=Q_{xx}(s,X(s)),\q\BC(s)=S(s),\qq
s\in[0,T].$$
Then by (\ref{5.15})--(\ref{5.16}), we can apply Proposition 5.2 to
obtain the uniform positive definiteness of $D_{uu}J(t,x;u(\cd))$.
\endpf

\ms

Recall that maps $x\mapsto Q(t,x)$ and $x\mapsto G(x)$ as {\it
semi-convex} maps if there is a constant $K>0$ such that
$$x\mapsto Q(t,x)+K|x|^2,\qq x\mapsto Q(x)+K|x|^2$$
are convex. It is clear that under (\ref{5.15}), $x\mapsto Q(t,x)$
and $x\mapsto G(x)$ are semi-convex. Hence, the associated problem
is referred to as a {\it linear semi-convex} problem. Our result
basically shows that the possible deviation from the convexity of
the maps $x\mapsto Q(t,x)$ and $x\mapsto G(x)$ could be possibly
compensated by the sufficient positive definiteness of $R(\cd)$.

\subsection{A more general case.}

We now impose the following conditions:
\bel{5.12}B(t,x)=B(t),\q
R(t,x)=R(t),\qq(t,x)\in[0,T]\times\dbR^n.\ee
Note that we still allow $x\mapsto(A(t,x),Q(t,x),S(t,x),G(x))$ to be
nonlinear. In the current case, we have
\bel{}\left\{\2n\ba{ll}
\ns\ds\BA(s)=A_x(s,X(s)),\\
\ns\ds\BA_1(s)=\sum_{i=1}^nY^i(s)A^i_{xx}(s,X(s))+Q_{xx}(s,X(s))
+\2n\sum_{j=1}^mu^j(s)S^j_{xx}(s,X(s)),\\
\ns\ds\BC(s)=S_x(s,X(s)).\ea\right.\ee
Also,
$$\ba{ll}
\ns\ds\int_t^T\lan\,[D_{uu}J(t,x;u(\cd))v(\cd)](s),v(s)\ran ds\\
\ns\ds\ge\2n\int_t^T\lan(1-\a)R(s)v(s),v(s)\ran ds-\lan\bar
G\int_t^T\F_\BA(T,r)B(r)v(r)dr,\int_t^T\F_\BA(T,r)B(r)v(r)dr\ran\\
\ns\ds~-\int_t^T\lan\bar
Q(s)\int_t^s\F_\BA(s,r)B(r)v(r)dr\int_t^s\F_\BA(s,r)B(r)v(r)dr\ran ds\\
\ns\ds~+\lan(G_{xx}(X(T))+\bar G)\int_t^T
\F_\BA(T,r)B(r)v(r)dr,\int_t^T\F_\BA(T,r)B(r)v(r)dr\ran\\
\ns\ds~+\2n\int_t^T\3n\1n\lan\1n\big[\BA_1(s)\1n+\1n\bar
Q(s)\1n-\1n\a^{\1n-1} S_x(s,X(s))^T\1n
R(s)^{\1n-1}S_x(s,X(s))\big]\2n\int_t^s\3n\F_\BA(s,r)B(r)v(r)dr,\2n\int_t^s\3n\F_\BA(s,r)B(r)v(r)dr
\1n\ran\\
\ns\ds=\int_t^T\lan\{(1-\a)R(s)-[\h G(t)+\h Q(s,t)]I\}v(s),v(s)\ran
ds+\dbI_2+\dbI_3,\ea$$
where $\h G(\cd)$ and $\h Q(\cd)$ are defined by the
following:
\bel{5.20}\left\{\2n\ba{ll}
\ns\ds\h G(t)=\[\int_t^T\int_t^T|B(s)^T\F_\BA(T,s)^T\bar
G\F_\BA(T,r)B(r)|^2drds\]^{1\over2},\\
\ns\ds\h
Q(s,t)=\int_s^T\[\int_t^\t\int_t^\t|B(r)^T\F_\BA(\t,r)^T\bar
Q(\t)\F_\BA(\t,r')B(r')|^2dr'dr\]^{1\over2}d\t,\ea\right.\ee
and
$$\left\{\2n\ba{ll}
\ns\ds\dbI_2=\lan(G_{xx}(X(T))+\bar G)\int_t^T
\F_\BA(T,r)B(r)v(r)dr,\int_t^T\F_\BA(T,r)B(r)v(r)dr\ran\\
\ns\ds\dbI_3\1n=\2n\int_t^T\3n\1n\lan\1n\big[\BA_1(s)\1n+\1n\bar
Q(s)\1n-\1n\a^{\1n-1}\1n S_x(s,X(s))^T\1n R(s)^{-1}\1n
S_x(s,X(s))\big]\2n\int_t^s\3n\F_\BA(s,r)B(r)v(r)dr,\1n
\int_t^s\3n\F_\BA(s,r)B(r)v(r)dr\1n\ran.\ea\right.$$
Naturally, we may still assume
$$G_{xx}(x)+\bar G\ge0,\qq\forall x\in\dbR^n.$$
To ensure $\dbI_3\ge0$, we need to take a closer look at the
involved terms. Note that in the current case, $\BA_1(s)$ involves
$(X(\cd),Y(\cd))$ and $u(\cd)$, unless $A_{xx}^i(s,x)=0$ and
$S^j_{xx}(s,x)=0$. Basically, we hope to get a uniform boundedness
from below. On the other hand, it is known that we need the positive
definiteness of $D_{uu}J(t,x;u(\cd))$ for any candidate $u(\cd)$ of
optimal control only, rather than any admissible control. Hence, we
restrict $u(\cd)$ as follows:
$$u(s)=-R(s)^{-1}\[B(s)^TY(s)+S(s,X(s))\],\qq s\in[t,T],$$
where $(X(\cd),Y(\cd))$ is a solution to the following:
\bel{5.13}\2n\3n\left\{\2n\ba{ll}
\ns\ds\dot X(s)=A(s,X(s))-B(s)R(s)^{-1}\[B(s)^TY(s)+S(s,X(s))\],\\
\ns\ds\dot
Y(s)=-\[A_x(s,X(s))-B(s)R(s)^{-1}S_x(s,X(s))\]^TY(s)\\
\ns\ds\qq\qq-\[Q_x(s,X(s))-S(s,X(s))^TR(s)^{-1}S_x(s,X(s))\]^T,\q s\in[t,T],\\
\ns\ds X(t)=x,\qq Y(T)=G_x(X(T))^T.\ea\right.\ee
Consequently, we have
$$\ba{ll}
\ns\ds\BA_1(s)=\sum_{i=1}^nY^i(s)\[A^i_{xx}(s,X(s))-\sum_{j=1}^me_j^TR(s)^{-1}\wt B^i(s)S_{xx}^i(s,X(s))\]+Q_{xx}(s,X(s))\\
\ns\ds\qq\qq-\sum_{j=1}^m\[e_j^TR(s)^{-1}
S(s,X(s))\]S^j_{xx}(s,X(s)).\ea$$
We now would like to explore the possibility of
\bel{5.22}\ba{ll}
\ns\ds0\le\BA_1(s)\1n+\1n\bar Q(s)\1n-\1n\a^{\1n-1}\1n
S_x(s,X(s))^T\1n R(s)^{-1}\1n S_x(s,X(s))\\
\ns\ds\q=Q_{xx}(s,X(s))+\bar
Q(s)-\a^{-1}S_x(s,X(s))^TR(s)^{-1}S_x(s,X(s))\\
\ns\ds\qq-\sum_{j=1}^m\[e_j^TR(s)^{-1}
S(s,X(s))\]S^j_{xx}(s,X(s))\\
\ns\ds\qq+\sum_{i=1}^nY^i(s)\[A^i_{xx}(s,X(s))
-\sum_{j=1}^me_j^TR(s)^{-1}\wt B^i(s)S_{xx}^i(s,X(s))\],\ea\ee
for some $\bar Q(\cd)$. If we are above to show that
\bel{Y<X}|Y(s)|\le K_0\big(1+|X(s)|\big),\qq s\in[t,T],\ee
then (\ref{5.22}) is guaranteed by the following:
\bel{}\ba{ll}
\ns\ds Q_{xx}(s,x)+\bar
Q(s)-\a^{-1}S_x(s,x)^TR(s)^{-1}S_x(s,x)-\sum_{j=1}^m\[e_j^TR(s)^{-1}
S(s,x)\]S^j_{xx}(s,x)\\
\ns\ds+\sum_{i=1}^nK_0\big(1+|x|\big)\[A^i_{xx}(s,x)
-\sum_{j=1}^me_j^TR(s)^{-1}\wt
B^i(s)S_{xx}^i(s,x)\]\ge0,\qq\forall(s,x)\in[0,T],\ea\ee
which is practically checkable. Interestingly, in the case that
$x\mapsto(A(s,x),S(s,x))$ is linear, the above is reduced to the
first condition in (\ref{5.15}), and (\ref{Y<X}) is not necessary.

\ms

Now, let us look at conditions under which (\ref{Y<X}) holds.

\ms

\bf Lemma 5.5. \sl Let {\rm (H1)--(H3)} and $(\ref{5.12})$ hold. Let
\bel{5.14}\left\{\2n\ba{ll}
\ns\ds|A_x(t,x)-B(t)R(t)^{-1}S_x(t,x)|\le\bar L_A,\qq
|B(t)R(t)^{-1}B(t)^T|\le\bar L_B,\\
\ns\ds|Q_x(t,x)-S(t,x)^TR(t)^{-1}S_x(t,x)|\le Q_0+\bar
L_Q|x|,\qq|G_x(x)|\le|G_x(0)|+\bar L_G|x|,\\
\ns\ds\qq\qq\qq\qq\qq\qq\qq\qq\qq\forall(t,x)\in[0,T]\times\dbR^n,\ea\right.\ee
for some constants $\bar L_A,\bar L_B,\bar L_Q,\bar L_G\ge0$.
Suppose
\bel{T()<1}T\big[\bar L_A+\bar L_B(\bar L_Q+\bar L_G)e^{\bar
L_AT}\big]<1.\ee
Then there exists an absolute constant $K_0>0$, independent of
$(t,x)$, such that for any $(t,x)\in[0,T)\times\dbR^n$, any solution
$(X(\cd),Y(\cd))$ of $(\ref{5.13})$ satisfies
\bel{}|Y(s)|\le K_0\(1+|X(s)|\),\qq s\in[t,T].\ee

\rm

\it Proof. \rm In what follows, we denote
$$\left\{\2n\ba{ll}
\ns\ds A_0=\max_{t\in[0,T]}|A(t,0)-B(t)R(t)^{-1}S(t,0)|,\\
\ns\ds
Q_0=\max_{t\in[0,T]}|Q_x(t,0)-S(t,0)^TR(t)^{-1}S_x(t,0)|.\ea\right.$$
For given $(t,x)\in[0,T)\times\dbR^n$, let $(X(\cd),Y(\cd))$ be a
solution to (\ref{5.13}). Let $\Psi(\cd\,,\cd)$ be the fundamental
matrix of $A_x(\cd\,,X(\cd))-B(\cd)R(\cd)^{-1}S_x(\cd\,,X(\cd))$,
i.e.,
$$\left\{\2n\ba{ll}
\ns\ds\dot\Psi(s,\t)=\[A_x(s,X(s))-B(s)R(s)^{-1}S_x(s,X(s))\]\Psi(s,\t),\qq
\t,s\in[t,T],\\
\ns\ds\Psi(\t,\t)=I,\ea\right.$$
Then by the first condition in (\ref{5.14}), we have
$$|\Psi(s,\t)|\le e^{\bar L_A(s-\t)},\qq\forall t\le\t\le s\le T.$$
Next, $Y(\cd)$ admits the following representation:
$$Y(s)\1n=\1n\Psi(T,s)^T\1n
G_x(X(T))^T\2n+\1n\int_s^T\2n\Psi(r,s)^T\1n\[Q_x(r,X(r))\1n-\1n
S(r,X(r))^T\1n R(r)^{-1}S_x(r,X(r))\]^T\1n dr.$$
Thus,
$$\ba{ll}
\ns\ds|Y(s)|\1n\le\1n e^{\bar L_A(T-s)}|G_x(0)|\1n+\1n{Q_0(e^{\bar
L_A(T-s)}\1n-\1n1)\over\bar L_A}\1n+\1n\bar L_Ge^{\bar
L_A(T-s)}|X(T)|\1n+\1n\bar L_Q\int_s^Te^{\bar L_A(r-s)}
|X(r)|dr\\
\ns\ds\qq~\le K_1+\bar L_Ge^{\bar L_AT}|X(T)|+\bar L_Qe^{\bar
L_AT}\int_s^T|X(r)|dr,\ea$$
where
$$K_1=e^{\bar L_AT}|G_x(0)|\1n+\1n{Q_0(e^{\bar L_AT}\1n-\1n1)\over\bar L_A}\,.$$
On the other hand, from
$$\left\{\2n\ba{ll}
\ns\ds\dot X(s)=A(s,X(s))-B(s)R(s)^{-1}\[B(s)^TY(s)+S(s,X(s))\],\q s\in[t,T],\\
\ns\ds X(t)=x,\ea\right.$$
for any $t\le s\le\t\le T$, we have
$$\ba{ll}
\ns\ds|X(\t)|\1n\le\1n|X(s)|\1n+\3n\int_s^\t\3n\(|A(r,X(r))\1n-\1n
B(r)R(r)^{-1}S(r,X(r))|\1n
+\1n|B(r)R(r)^{-1}B(r)^T|\,|Y(r)|\)dr\\
\ns\ds\le|X(s)|\1n+\2n\int_s^\t\1n\[A_0\1n+\1n\bar
L_A|X(r)|\1n+\1n\bar L_B\(K_1\1n+\1n\bar L_Ge^{\bar
L_AT}|X(T)|\1n+\1n\bar L_Qe^{\bar L_AT}
\int_r^T|X(r')|dr'\)\]dr\\
\ns\ds\le|X(s)|+(A_0+\bar L_BK_1)(\t-s)+\bar
L_A\int_s^\t|X(r)|dr+\bar L_B\bar L_G
e^{\bar L_AT}(\t-s)|X(T)|\\
\ns\ds\qq\qq+\bar L_B\bar L_Q\int_s^\t
e^{\bar L_AT}\int_r^T|X(r')|dr'dr\\
\ns\ds\le(A_0+\bar L_BK_1)T+|X(s)|+(\bar L_A+\bar L_B\bar
L_QTe^{\bar L_AT})\int_s^T|X(r)|dr+\bar L_B\bar L_GTe^{\bar
L_AT}|X(T)|.\ea$$
Then with $\t=T$, we have
$$\ba{ll}
\ns\ds|X(T)|\le(A_0+\bar L_BK_1)T+|X(s)|+(\bar L_A+\bar L_B\bar
L_QTe^{\bar L_AT})\int_s^T|X(r)|dr+\bar L_B\bar L_GTe^{\bar
L_AT}|X(T)|.\ea$$
Hence, under condition (\ref{T()<1}), one has
$$|X(T)|\le{(A_0+\bar L_BK_1)T\over1-\bar L_B\bar
L_GTe^{\bar L_AT}}+{1\over1-\bar L_B\bar L_GTe^{\bar
L_AT}}|X(s)|+{\bar L_A+\bar L_B\bar L_QTe^{\bar L_AT}\over1-\bar
L_B\bar L_GTe^{\bar L_AT}}\int_s^T|X(r)|dr.$$
Consequently,
$$\ba{ll}
\ns\ds\int_s^T\3n|X(r)|dr\1n\le\1n(A_0\1n+\1n\bar L_BK_1)T^2\1n+\1n
T|X(s)|\1n+\1n(\bar L_A\1n+\1n\bar L_B\bar L_QTe^{\bar
L_AT})T\2n\int_s^T\3n|X(r)|dr\1n+\1n\bar
L_B\bar L_GT^2e^{\bar L_AT}|X(T)|\\
\ns\ds\le(A_0+\bar L_BK_1)T^2+T|X(s)|+(\bar L_A+\bar L_B\bar
L_QTe^{\bar L_AT})T\int_s^T|X(r)|dr\\
\ns\ds\q+\bar L_B\bar L_GT^2e^{\bar L_AT}\[{(A_0+\bar
L_BK_1)T\over{1-\bar L_B\bar L_GTe^{\bar L_AT}}}+{1\over1-\bar
L_B\bar L_GTe^{\bar L_AT}}|X(s)|+{\bar L_A\1n+\1n\bar L_B\bar
L_QTe^{\bar L_AT}\over1-\bar L_B\bar
L_GTe^{\bar L_AT}}\1n\int_s^T\3n|X(r)|dr\]\\
\ns\ds=(A_0+\bar L_BK_1)T^2\[1+{\bar L_B\bar L_GTe^{\bar
L_AT}\over1-\bar L_B\bar L_GTe^{\bar L_AT}}\]+T\[1+{\bar L_B\bar
L_GTe^{\bar L_AT}\over1-
\bar L_B\bar L_GTe^{\bar L_AT}}\]|X(s)|\\
\ns\ds\q+(\bar L_A+\bar L_B\bar L_QTe^{\bar L_AT})T\[1+{\bar L_B\bar
L_GTe^{\bar L_AT}\over
1-\bar L_B\bar L_GTe^{\bar L_AT}}\]\int_s^T|X(r)|dr\\
\ns\ds={(A_0+\bar L_BK_1)T^2\over1-\bar L_B\bar L_GTe^{\bar
L_AT}}+{T\over1-\bar L_B\bar L_GTe^{\bar L_AT}}|X(s)|+{(\bar
L_A+\bar L_B\bar L_Qe^{\bar L_AT})T\over1-\bar L_B\bar L_GTe^{\bar
L_AT}}\int_s^T|X(r)|dr.\ea$$
Therefore, under condition (\ref{T()<1}), one has
$$\ba{ll}
\ns\ds\int_s^T|X(r)|dr\le{(A_0+\bar
L_BK_1)T^2+T|X(s)|\over1-T\big[\bar L_A+\bar L_B(\bar L_Q+\bar
L_G)e^{\bar L_AT}\big]}\equiv K_2(A_0+\bar L_BK_1)T+K_2|X(s)|,\ea$$
where
$$K_2={T\over1-T\big[\bar L_A+\bar L_B(\bar L_Q+\bar L_G)e^{\bar L_AT}\big]}.$$
Consequently,
$$\ba{ll}
\ns\ds|X(T)|\le{(A_0+\bar L_BK_1)T\over{1-\bar L_B\bar L_GTe^{\bar
L_AT}}}+{1\over1-\bar L_B\bar L_GTe^{\bar L_AT}}|X(s)|+{\bar
L_A+\bar L_B\bar L_QTe^{\bar L_AT}\over1-\bar L_B\bar
L_GTe^{\bar L_AT}}\int_s^T|X(r)|dr\\
\ns\ds\le{(A_0+\bar L_BK_1)T\over{1-\bar L_B\bar
L_GTe^{\bar L_AT}}}+{1\over1-\bar L_B\bar L_GTe^{\bar L_AT}}|X(s)|\\
\ns\ds\qq+{\bar L_A+\bar L_B\bar L_QTe^{\bar L_AT}\over1-\bar
L_B\bar L_GTe^{\bar L_AT}}
\[K_2(A_0+\bar L_BK_1)T+K_2|X(s)|\]\\
\ns\ds={(A_0+\bar L_BK_1)T\over{1-\bar L_B\bar L_GTe^{\bar
L_AT}}}\[1+K_2(\bar L_A+\bar L_B\bar L_QTe^{\bar L_AT})\]
+{1+K_2(\bar L_A+\bar L_B\bar L_QTe^{\bar L_AT})
\over1-\bar L_B\bar L_GTe^{\bar L_AT}}|X(s)|\\
\ns\ds=K_3(A_0+\bar L_BK_1)T+K_3|X(s)|,\ea$$
with
$$K_3={1+K_2(\bar L_A+\bar L_B\bar L_QTe^{\bar L_AT})
\over1-\bar L_B\bar L_GTe^{\bar L_AT}}.$$
Hence,
$$\ba{ll}
\ns\ds|Y(s)|\le K_1+\bar L_Ge^{\bar L_AT}|X(T)|+\bar
L_Qe^{\bar L_AT}\int_s^T|X(r)|dr\\
\ns\ds\le K_1+\bar L_Ge^{\bar L_AT}\[K_3(A_0+\bar
L_BK_1)T+K_3|X(s)|\]+\bar L_Qe^{\bar L_AT}\[K_2(A_0+\bar
L_BK_1)T+K_2|X(s)|\]\\
\ns\ds=K_1+e^{\bar L_AT}(A_0+\bar L_BK_1)T(\bar L_GK_3+\bar L_QK_2)
+e^{\bar L_AT}(\bar L_GK_3+\bar L_QK_2)|X(s)|\equiv
K_4+K_5|X(s)|,\ea$$
with
$$\left\{\2n\ba{ll}
\ns\ds K_4=K_1+(\bar L_GK_3+\bar L_QK_2)(A_0+\bar L_BK_1)Te^{\bar L_AT},\\
\ns\ds K_5=(\bar L_GK_3+\bar L_QK_2)e^{\bar L_AT}.\ea\right.$$
This proves our lemma with $K_0=K_4\vee K_5$.

\ms

Condition (\ref{T()<1}) tells us that (\ref{Y<X}) holds if $T>0$ is
not too large, and the found constant $K_0$ depends on all the
constants $\bar L_A,\bar L_B,\bar L_Q,\bar L_G,A_0,Q_0$, as well as
the time duration $T$. The following result, under some different
conditions, shows that sometimes, $T>0$ could be arbitrarily large.

\ms

\bf Lemma 5.6. \sl Let {\rm (H1)--(H3)} and $(\ref{5.12})$ hold. Let
\bel{5.28}\left\{\2n\ba{ll}
\ns\ds\big[A_x(t,x)\1n-\1n B(t)R(t)^{-1}S_x(t,x)\big]
\1n+\1n\big[A_x(t,x)\1n-\1n B(t)R(t)^{-1}S_x(t,x)\big]^T\3n\le\2n
-2L_0I,\\
\ns\ds|A_x(t,x)-B(t)R(t)^{-1}S_x(t,x)|\le\bar L_A,\qq
|B(t)R(t)^{-1}B(t)^T|\le\bar L_B,\\
\ns\ds|Q_x(t,x)-S(t,x)^TR(t)^{-1}S_x(t,x)|\le\bar
L_Q|x|,\qq|G_x(x)|\le\bar
L_G|x|,\ea\right.\q\forall(t,x)\in[0,T]\times\dbR^n,\ee
for some constants $\bar L_A,\bar L_B,\bar L_G,\bar L_Q\ge0$ and
$L_0>0$. Suppose
\bel{A=S=0}A(t,0)=S(t,0)=0,\qq\forall t\in[0,T],\ee
and
\bel{<}\left\{\ba{ll}
\ns\ds2L_0^2-\bar L_B^2\bar L_G^2>0,\qq L_0^4-\bar L_B^2\bar
L_Q^2>0,\\
\ns\ds2(2L_0^2-\bar L_B^2\bar L_G^2)(L_0^4-\bar L_B^2\bar
L_Q^2)-\bar L_B^4\bar L_Q^2\bar L_G^2>0.\ea\right.\ee
For any $(t,x)\in[0,T)\times\dbR^n$, let $(X(\cd),Y(\cd))$ be a
solution of $(\ref{5.13})$. Then there exists an absolute constant
$K_0>0$, independent of $(t,x)$ and $T$, such that
\bel{5.31}|Y(s)|\le K_0|X(s)|,\qq s\in[t,T].\ee

\rm

\it Proof. \rm For given $(t,x)\in[0,T)\times\dbR^n$, let
$(X(\cd),Y(\cd))$ be a solution to (\ref{5.13}). Then making use of
the first condition in (\ref{5.28}),
$$\ba{ll}
\ns\ds|Y(s)|^2=|G_x(X(T))|^2+2\int_s^T\(\lan
Y(r),\big[A_x(r,X(r))-B(r)R(r)^{-1}S_x(r,X(r))\big]^TY(r)\ran\\
\ns\ds\qq\qq\qq+\lan
Y(r),\big[Q_x(r,X(r))-S(r,X(r))^TR(r)^{-1}S_x(r,X(r))\big]^T\ran\)dr\\
\ns\ds\qq\q\le\bar L_G^2|X(T)|^2-2L_0\int_s^T|Y(r)|^2dr+2\bar
L_Q\int_s^T|Y(r)||X(r)|dr\\
\ns\ds\qq\q\le\bar L_G^2|X(T)|^2-L_0\int_s^T|Y(r)|^2+\int_s^T{\bar
L_Q^2\over L_0}|X(r)|^2dr.\ea$$
Hence, by Gronwall's inequality,
$$|Y(s)|^2\le\bar L_G^2e^{-L_0(T-s)}|X(T)|^2+{\bar L_Q^2\over L_0}
\int_s^Te^{-L_0(r-s)}|X(r)|^2dr.$$
Next, by the first condition in (\ref{5.28}) again, together with
(\ref{A=S=0}), we have
$$\ba{ll}
\ns\ds\lan
x,A(t,x)-B(t)R(t)^{-1}S(t,x)\ran=\lan\(\int_0^1\big[A_x(t,\b
x)-B(t)R(t)^{-1}S_x(t,\b x)\big]d\b\)x,x\ran\le-L_0|x|^2.\ea$$
Thus, from
$$\left\{\2n\ba{ll}
\ns\ds\dot X(s)=A(s,X(s))-B(s)R(s)^{-1}\[B(s)^TY(s)+S(s,X(s))\],\q s\in[t,T],\\
\ns\ds X(t)=x,\ea\right.$$
for any $t\le s\le\t\le T$, we have
$$\ba{ll}
\ns\ds|X(\t)|^2=|X(s)|^2+2\int_s^\t\(\lan
X(r),A(r,X(r))-B(r)R(r)^{-1}S(r,X(r))\ran\\
\ns\ds\qq\qq\qq\qq-\lan X(r),B(r)R(r)^{-1}B(r)^TY(r)\ran\)dr\\
\ns\ds\le|X(s)|^2-2\int_s^\t\(L_0|X(r)|^2-\bar
L_B|X(r)||Y(r)|\)dr\\
\ns\ds\le|X(s)|^2+\int_s^\t\(-L_0|X(r)|^2+{\bar L_B^2\over
L_0}|Y(r)|^2\)dr.\ea$$
Then, by Gronwall's inequality,
\bel{5.32}\ba{ll}
\ns\ds|X(\t)|^2\le e^{-L_0(\t-s)}|X(s)|^2+{\bar L_B^2\over
L_0}\int_s^\t e^{-L_0(\t-r)}|Y(r)|^2dr\\
\ns\ds\le e^{-L_0(\t-s)}|X(s)|^2+{\bar L_B^2\over L_0}\int_s^\t
e^{-L_0(\t-r)}\(\bar L_G^2e^{-L_0(T-r)}|X(T)|^2+{\bar L_Q^2\over
L_0}\int_r^Te^{-L_0(r'-r)}|X(r')|^2dr'\)dr\\
\ns\ds=e^{-L_0(\t-s)}|X(s)|^2+{\bar L_B^2\bar
L_G^2(e^{-L_0(T-\t)}-e^{-L_0(T+\t-2s)})\over
2L_0^2}|X(T)|^2\\
\ns\ds\qq\qq+{\bar L_B^2\bar L_Q^2\over L_0^2}\int_s^\t
e^{-L_0(\t-r)}\(\int_r^Te^{-L_0(r'-r)}|X(r')|^2dr'\)dr\\
\ns\ds=e^{-L_0(\t-s)}|X(s)|^2+{\bar L_B^2\bar
L_G^2(e^{-L_0(T-\t)}-e^{-L_0(T+\t-2s)})\over
2L_0^2}|X(T)|^2\\
\ns\ds\qq\qq+{\bar L_B^2\bar L_Q^2\over2L_0^3}\int_s^T
(e^{-L_0[\t+r-2(\t\land r)]}-e^{-L_0(\t+r-2s)})|X(r)|^2dr\\
\ns\ds\le e^{-L_0(\t-s)}|X(s)|^2+{\bar L_B^2\bar
L_G^2e^{-L_0(T-\t)}\over 2L_0^2}|X(T)|^2+{\bar L_B^2\bar
L_Q^2\over2L_0^3}\int_s^T e^{-L_0[\t+r-2(\t\land
r)]}|X(r)|^2dr.\ea\ee
Integrating the above over $[s,T]$, we obtain
$$\ba{ll}
\ns\ds\int_s^T|X(\t)|^2d\t\le\int_s^T\[e^{-L_0(\t-s)}|X(s)|^2+{\bar
L_B^2\bar L_G^2e^{-L_0(T-\t)}\over 2L_0^2}|X(T)|^2\\
\ns\ds\qq\qq\qq\qq\qq+{\bar L_B^2\bar
L_Q^2\over2L_0^3}\int_s^T e^{-L_0[\t+r-2(\t\land r)]}|X(r)|^2dr\]d\t\\
\ns\ds\le{1-e^{L_0(T-s)}\over L_0}|X(s)|^2+{\bar L_B^2\bar
L_G^2(1-e^{-L_0(T-s)})\over2L_0^3}|X(T)|^2\\
\ns\ds\qq+{\bar L_B^2\bar L_Q^2\over2L_0^3}\int_s^T\(\int_s^\t
e^{-L_0(\t-r)}|X(r)|^2dr+\int_\t^Te^{-L_0(r-\t)}|X(r)|^2dr\)d\t\\
\ns\ds\le{1\over L_0}|X(s)|^2+{\bar L_B^2\bar
L_G^2\over2L_0^3}|X(T)|^2+{\bar L_B^2\bar
L_Q^2\over2L_0^3}\[\int_s^T\3n\(\int_r^Te^{-L_0(\t-r)}d\t\)|X(r)|^2dr\\
\ns\ds\qq\qq\qq\qq\qq\qq\qq\qq\qq+\int_s^T\(\int_s^re^{-L_0(r-\t)}d\t\)|X(r)|^2dr\]\\
\ns\ds={1\over L_0}|X(s)|^2+{\bar L_B^2\bar
L_G^2\over2L_0^3}|X(T)|^2+{\bar L_B^2\bar
L_Q^2\over2L_0^4}\[\int_s^T\((1-e^{-L_0(T-r)}+(1-e^{-L_0(r-s)})\)
|X(r)|^2dr\]\\
\ns\ds\le{1\over L_0}|X(s)|^2+{\bar L_B^2\bar
L_G^2\over2L_0^3}|X(T)|^2+{\bar L_B^2\bar L_Q^2\over L_0^4}\int_s^T
|X(r)|^2dr.\ea$$
Thus, under condition (\ref{<}), one has
$$\int_s^T|X(r)|^2dr\le{L_0^3\over L_0^4-\bar L_B^2\bar
B_Q^2}|X(s)|^2+{L_0\bar L_B^2\bar L_G^2\over2(L_0^4-\bar L_B^2\bar
L_Q^2)}|X(T)|^2.$$
On the other hand, taking $\t=T$ in (\ref{5.32}), we get
$$|X(T)|^2\le e^{-L_0(T-s)}|X(s)|^2+{\bar L_B^2\bar
L_G^2\over 2L_0^2}|X(T)|^2+{\bar L_B^2\bar L_Q^2\over2L_0^3}\int_s^T
e^{-L_0(T-r)}|X(r)|^2dr.$$
Thus, by (\ref{<}),
$$\ba{ll}
\ns\ds|X(T)|^2\le{2L_0^2e^{-L_0(T-s)}\over2L_0^2-\bar L_B^2\bar
L_G^2}|X(s)|^2+{\bar L_B^2\bar L_Q^2\over L_0(2L_0^2-\bar L_B^2\bar
L_G^2)}\int_s^Te^{-L_0(T-r)}|X(r)|^2dr\\
\ns\ds\le{2L_0^2\over2L_0^2-\bar L_B^2\bar L_G^2}|X(s)|^2+{\bar
L_B^2\bar L_Q^2\over L_0(2L_0^2-\bar L_B^2\bar
L_G^2)}\int_s^T|X(r)|^2dr\\
\ns\ds\le{2L_0^2\over2L_0^2-\bar L_B^2\bar L_G^2}|X(s)|^2+{\bar
L_B^2\bar L_Q^2\over L_0(2L_0^2-\bar L_B^2\bar L_G^2)}\[{L_0^3\over
L_0^4-\bar L_B^2\bar B_Q^2}|X(s)|^2+{L_0\bar
L_B^2\bar L_G^2\over2(L_0^4-\bar L_B^2\bar L_Q^2)}|X(T)|^2\]\\
\ns\ds=\[{2L_0^2\over2L_0^2-\bar L_B^2\bar L_G^2} +{L_0^2\bar
L_B^2\bar L_Q^2\over(2L_0^2-\bar L_B^2\bar L_G^2)(L_0^4-\bar
L_B^2\bar L_Q^2)}\]|X(s)|^2+{\bar L_B^4\bar L_Q^2\bar
L_G^2\over2(L_0^4-\bar L_B^2\bar L_Q^2)(2L_0^2-\bar L_B^2\bar
L_G^2)}|X(T)|^2.\ea$$
Hence, by (\ref{<}), we obtain
$$|X(T)|^2\le{4L_0^2(L_0^4-\bar L_B^2\bar L_Q^2)+2L_0^2\bar
L_B^2\bar L_Q^2\over2(L_0^4-\bar L_B^2\bar L_Q^2)(2L_0^2-\bar
L_B^2\bar L_G^2)-\bar L_B^4\bar L_Q^2\bar L_G^2}|X(s)|^2,$$
and
$$\int_s^T|X(r)|^2dr\le\[{L_0^3\over L_0^4-\bar L_B^2\bar
B_Q^2}+{4L_0^2(L_0^4-\bar L_B^2\bar L_Q^2)+2L_0^2\bar L_B^2\bar
L_Q^2\over2(L_0^4-\bar L_B^2\bar L_Q^2)(2L_0^2-\bar L_B^2\bar
L_G^2)-\bar L_B^4\bar L_Q^2\bar L_G^2}\]|X(s)|^2.$$
Combining the above, we finally get
$$\ba{ll}
\ns\ds|Y(s)|^2\le\bar L_G^2|X(T)|^2+{\bar L_Q^2\over L_0}
\int_s^T|X(r)|^2dr\\
\ns\ds\le\[{\bar L_G^2\big[4L_0^2(L_0^4-\bar L_B^2\bar
L_Q^2)+2L_0^2\bar L_B^2\bar L_Q^2\big]\over2(L_0^4-\bar L_B^2\bar
L_Q^2)(2L_0^2-\bar L_B^2\bar L_G^2)-\bar L_B^4\bar L_Q^2\bar
L_G^2}\\
\ns\ds\qq+{L_0^2\bar L_Q^2\over L_0^4-\bar L_B^2\bar
L_Q^2}+{4L_0\bar L_Q^2(L_0^4-\bar L_B^2\bar L_Q^2)+2L_0\bar
L_B^2\bar L_Q^4\over2(L_0^4-\bar L_B^2\bar L_Q^2)(2L_0^2-\bar
L_B^2\bar L_G^2)-\bar L_B^4\bar L_Q^2\bar L_G^2}\)\]|X(s)|^2\\
\ns\ds=\[{L_0^2\bar L_Q^2\over L_0^4-\bar L_B^2\bar L_Q^2}+{(\bar
L_G^2+L_0\bar L_Q^2)\big[4L_0^2(L_0^4-\bar L_B^2\bar
L_Q^2)+2L_0^2\bar L_B^2\bar L_Q^2\big]\over2(L_0^4-\bar L_B^2\bar
L_Q^2)(2L_0^2-\bar L_B^2\bar L_G^2)-\bar L_B^4\bar L_Q^2\bar
L_G^2}\]|X(s)|^2.\ea$$
Then our conclusion follows. \endpf

\ms

If (\ref{A=S=0}) is not assumed and the last line in (\ref{5.28}) is
replaced by the second line of (\ref{5.14}), then, with more
complicated-looking estimates, we will have
$$|Y(s)|\le K_0(1+|X(s)|),\qq s\in[t,T],$$
instead of (\ref{5.31}), with $K_0$ also independent of $T$. For the
simplicity of presentation, we prefer not to give the details here.
Having the above two lemmas, we may state the following result whose
proof is clear.

\ms

\bf Theorem 5.7. \sl Let conditions of Lemma 5.5 or Lemma 5.6 hold.
Let $\bar G\in\dbS^n$ and $\bar Q:[0,T]\to\dbS^n$ such that for some
$\a\in(0,1)$ and $\d>0$,
\bel{}G_{xx}(x)+\bar G\ge0,\qq x\in\dbR^n,\ee
\bel{}R(s)-[\h G(t)+\h Q(s)]I\ge\d I,\qq 0\le t\le s\le T,\ee
with $\h G(\cd)$ and $\h Q(\cd)$ defined as in $(\ref{5.9})$, and
\bel{5.20}\ba{ll}
\ns\ds Q_{xx}(s,x)+\bar
Q(s)-K_0(1+|x|)\Big|\sum_{i=1}^nA^i_{xx}(s,x)
-\sum_{j=1}^me_j^TR(s)^{-1}\wt B^i(s)^TS^j_{xx}(s,x)\Big|I\\
\ns\ds\qq-\sum_{j=1}^m\[e_j^T
R(s)^{-1}S(s,x)\]S^j_{xx}(s,x)-(1-\a)^{-1}S_x(s,x)^TR(s)^{-1}S_x(s,x)\ge0,\ea\ee
with $K_0>0$ being obtained in Lemma 5.5 or Lemma 5.6.

\ms

\rm

From (\ref{5.20}), we see that due to the nonlinearity of $x\mapsto
(A(t,x),S(t,x))$, we basically need the semi-convexity of $x\mapsto
Q(t,x)$ and the sufficient positive definiteness of $R(\cd)$ to
compensate.

\section{Quasi-Riccati Equation.}

Let us keep condition (\ref{5.12}). We have seen that under certain
conditions $D_{uu}J(t,x;u^*(\cd))$ admits a bounded inverse at any
optimal control $u^*(\cd)\equiv u^*(\cd\,;t,x)$. When this is the
case, we have the following: Denoting $x_0=t$,
$$V(x_0,x)=J(x_0,x;u^*(\cd\,;x_0,x)),\qq D_uJ(x_0,x;u^*(\cd\,;x_0,x))=0.$$
Thus, for $0\le i\le n$,
$$u^*_{x_i}(\cd\,;x_0,x)=-D_{uu}J(x_0,x;u^*(\cd\,;x_0,x))^{-1}
D_uJ_{x_i}(x_0,x;u^*(\cd\,;x_0,x)).$$
Hence,
$$V_{x_i}(x_0,x)=J_{x_i}(x_0,x;u^*(\cd\,;x_0,x))+D_uJ(x_0,x;u^*(\cd\,;x_0,x))
u^*_{x_i}(\cd\,;x_0,x) =J_{x_i}(x_0,x;u^*(\cd\,;x_0,x)).$$
Further, for $0\le i,j\le n$,
$$\ba{ll}
\ns\ds V_{x_ix_j}(x_0,x)=J_{x_ix_j}(x_0,x;u^*(\cd\,;x_0,x))
+D_uJ_{x_i}(x_0,x;u^*(\cd\,;x_0,x))u^*_{x_j}(\cd\,;x_0,x)\\
\ns\ds\qq\qq\qq\qq+[D_{uu}J(x_0,x;u^*(\cd\,;x_0,x))u_{x_i}^*(\cd\,;x_0,x)]u_{x_j}^*
(\cd\,;x_0,x)\\
\ns\ds\qq\qq\q=J_{x_ix_j}(x_0,x;u^*(\cd\,;x_0,x)).\ea$$
Therefore, $V(\cd\,,\cd)$ is actually twice continuously
differentiable. Consequently, $V(\cd\,,\cd)$ satisfies the HJB
equation in the classical sense, and by the smoothness of the
coefficients, we can differentiate the equation once. Note that in
the current case, our HJB equation reads:
$$\left\{\2n\ba{ll}
\ns\ds V_t(t,x)+V_x(t,x)A(t,x)+Q(t,x)\\
\ns\ds\q-{1\over2}[V_x(t,x)B(t)\1n+\1n
S(t,x)^T]R(t)^{-1}[B(t,x)^T\1n V_x(t,x)^T\3n+\1n
S(t,x)]\1n=\1n0,\q(t,x)\1n\in\1n[0,T)\1n\times\1n
\dbR^n,\\
\ns\ds V(T,x)=G(x),\qq x\in\dbR^n.\ea\right.$$
Now, we define
$$P(t,x)=V_x(t,x)^T.$$
Then
$$P_x(t,x)=V_{xx}(t,x)=P_x(t,x)^T,\qq\forall(t,x)\in[0,T)\times\dbR^n,$$
and the following holds:
\bel{Quasi-R}\left\{\2n\ba{ll}
\ns\ds
P_t(t,x)+P_x(t,x)A(t,x)+A_x(t,x)^TP(t,x)+Q_x(t,x)^T\\
\ns\ds\q-[P_x(t,x)B(t)+S_x(t,x)^T]R(t)^{-1}
[B(t)^TP(t,x)+S(t,x)]=0,\q(t,x)\in[0,T]\times\dbR^n,\\
\ns\ds P(T,x)=G_x(x),\qq x\in\dbR^n.\ea\right.\ee
The above is called a {\it Quasi-Riccati equation} of Problem (AQ).
This is an extension of that presented in \cite{You 1997} for
linear-convex problems. We now have the following result.

\ms

\bf Proposition 6.1. \sl Let {\rm(H1)--(H3)} hold. Let
$(X^*(\cd),u^*(\cd))$ be an optimal pair of Problem (AQ). Suppose
the value function $V(\cd\,,\cd)$ of Problem (AQ) is twice
differentiable. Then $P(\cd\,,\cd)\equiv V_x(\cd\,,\cd)^T$ is a
solution to the quasi-Riccati equation $(\ref{Quasi-R})$, and the
optimal control $u^*(\cd)$ admits the following state feedback
representation:
\bel{}u^*(s)=-R(s)^{-1}\[B(s)^TP(s,X^*(s))+S(x,X^*(s))\],\qq
s\in[t,T].\ee

\ms

\it Proof. \rm It is known that if $(X^*(\cd),u^*(\cd))$ is an
optimal pair of Problem (AQ) for the initial pair
$(t,x)\in[0,T)\times\dbR^n$, and $Y(\cd)$ is the solution to the
corresponding adjoint equation, then
$$Y(s)=V_x(s.X^*(s))^T=P(s,X^*(s)),\qq s\in[t,T],$$
and
\bel{}\ba{ll}
\ns\ds
u^*(s)=-R(s)^{-1}\[B(s)^TY(s)+S(x,X^*(s))\]\\
\ns\ds\qq=-R(s)^{-1}\[B(s)^TP(s,X^*(s))+S(x,X^*(s))\],\qq
s\in[t,T].\ea\ee
This proves our conclusion. \endpf

\ms

In the case
$$A(t,x)=A(t)x,\q Q(t,x)={1\over2}\lan Q(t)x,x\ran,\q
S(t,x)=S(t)x,\q G(x)={1\over2}\lan Gx,x\ran$$
we see that
$$P(t,x)=P(t)x,\qq(t,x)\in[0,T]\times\dbR^n,$$
with $P(\cd)$ being the solution to the following:
$$\left\{\2n\ba{ll}
\ns\ds\dot P(t)\1n+\1n P(t)A(t)\1n+\1n A(t)^T\1n P(t)\1n+\1n
Q(t)\1n-\1n[P(t)B(t)\1n+\1n S(t)^T]R(t)^{-1}[B(t)^T\1n P(t)\1n+\1n
S(t)]\1n=\1n0,\q
t\in[0,T],\\
\ns\ds P(T)=G,\ea\right.$$
which is the Riccati equation for a standard LQ problem.

\ms

To conclude this section, we present two illustrative examples.

\ms

\bf Example 6.2. \rm Consider the following one-dimensional linear
controlled system:
$$\left\{\2n\ba{ll}
\ns\ds\dot X(s)=X(s)+u(s),\qq s\in[t,T],\\
\ns\ds X(t)=x,\ea\right.$$
with cost functional:
$$J(t,x;u(\cd))=\int_t^T\(-{1\over2}\cos^2X(s)+{1\over2}\rho X(s)^2\)ds-{1\over2}\sin^2X(T),$$
where $\rho>0$. This is a linear-semi-convex problem. According to
Proposition 5.4, we may choose
$$\bar Q(t)\equiv\bar G=1.$$
Then
$$\h G(t)=\int_t^Te^{2(T-s)}ds={e^{2(T-t)}-1\over2}\le\h G(0)={e^{2T}-1\over2},$$
and
$$\ba{ll}
\ns\ds\h Q(s)=\int_s^T\int_s^\t e^{2(\t-r)}drd\t={1\over2}\int_s^T\(e^{2(\t-s)}-1\)ds\\
\ns\ds\qq={e^{2(T-s)}-1\over4}-{T-s\over2}\le\h
Q(0)={e^{2T}-1\over4}-{T\over2}.\ea$$
Hence, under condition
$$\rho>{e^{2T}-1\over4}-{T\over2}+{e^{2T}-1\over2}={3(e^{2T}-1)\over4}-{T\over2},$$
we have the strict convexity of $u(\cd)\mapsto J(t,x;u(\cd))$.
Therefore, optimal control unique exists and the value function is
differentiable. In this case the optimal control $u^*(\cd)$ admits a
state feedback representation:
$$u^*(s)=-\rho^{-1}P(s,X^*(s)),\qq s\in[t,T],$$
with $P(\cd\,,\cd)$ solves the following quasi-Riccati equation:
$$\left\{\2n\ba{ll}
\ns\ds P_t(t,x)+xP_x(t,x)+P(t,x)+\sin2x-\rho^{-1}P_x(t,x)P(t,x)=0,\qq t\in[0,T],\\
\ns\ds P(T,x)=-\sin2x.\ea\right.$$

\ms

\bf Example 6.3. \rm Consider a one-dimensional controlled affine
system
$$\left\{\2n\ba{ll}
\ns\ds\dot X(s)=\sqrt{1+|X(s)|^2}+u(s),\qq s\in[t,T],\\
\ns\ds X(t)=x,\ea\right.$$
with cost functional
$$J(t,x;u(\cd))=\int_t^T\({Q\over2}|X(s)|^2+{R\over2}|u(s)|^2\)ds.$$
In this case, we may take
$$A_0=1,\q Q_0=0,\q\bar L_A=1,\q\bar L_B={1\over R},\q\bar L_G=0,\q\bar L_Q=q.$$
Then the first condition in (\ref{5.16}) automatically holds, and
the second condition reads
$${2R+Q(e^{2T}-1)\over2R}<{1\over T}.$$
This will be true if $T<1$  is small. Next, by looking at the proof
of Lemma 5.5, we see that
$$\ba{ll}
\ns\ds K_1=0,\qq K_2={RT\over2R(1-T)-QT(e^{2T}-1)},\\
\ns\ds
K_3={2R(1-T)-QT(e^{2T}-1)+RT\over2R(1-T)-QT(e^{2T}-1)},\\
\ns\ds K_4={QRT^2e^T\over2R(1-T)-QT(e^{2T}-1)},\\
\ns\ds K_5={QRTe^T\over2R(1-T)-QT(e^{2T}-1)}.\ea$$
Then
$$|Y(s)|\le K_4+K_5|X(s)|={QRTe^T(T+|X(s)|)\over2R(1-T)-QT(e^{2T}-1)}.$$
Now, if we assume
$$Q-{QRTe^T\over2R(1-T)-QT(e^{2T}-1)}{T+|x|\over(1+|x|^2)^{3\over2}}\ge0,\qq\forall x\in\dbR,$$
which can be achieved if $T>0$ is small, then Theorem 5.7 applies.
In the current case, the quasi-Riccati equation reads:
\bel{}\left\{\2n\ba{ll}
\ns\ds
P_t(t,x)\1n+\1n\sqrt{1\1n+\1n|x|^2}P_x(t,x)\1n+\1n{x\over\sqrt{1\1n+\1n|x|^2}}\,
P(t,x)\1n+\1n Q\1n-\1n P(t,x)P_x(t,x)\1n=0,
\q(t,x)\1n\in\1n[0,T]\1n\times\1n\dbR^n,\\
\ns\ds P(T,x)=0,\qq x\in\dbR^n.\ea\right.\ee
According to our result, under certain conditions (involving the
constant $K_0$), the above quasi-Riccati equation admits a solution
via which an optimal control admits a state feedback representation.

\section{Concluding Remarks.}

We have presented some very primitive results concerning what we
call the affine-quadratic optimal control problems, which are a
natural generalization of classical LQ problems, and also contains
linear-convex problems and linear-semi-convex problems. Our results
for linear state equation cover and substantially extend the known
results for LQ problems and linear-convex problems. Further, we have
some results for affine state equations. However, we see that there
are a lot challenging problems left open. Here are a couple of
these:

\ms

(i) Under our conditions, optimal controls exist and the optimality
system which is a two-point boundary value problem is always
solvable. It is a natural question if the corresponding
quasi-Riccati equation is always solvable? A technical question
relevant to this problem is: When the two-point boundary value
problem is always solvable, can one obtain an estimate
$$|Y(s)|\le K(1+|X(s)|),\qq s\in[t,T]$$
without additional restrictive conditions?

\ms

(ii) What happens if the dependence of $B(t,x)$ and $R(t,x)$ on $x$
is allowed? For such a situation, some new techniques might need to
be developed.

\ms

We expect to report some further relevant results in our future
publications.

\end{document}